\def\draft{n}
\documentclass{amsart}
\usepackage[headings]{fullpage}
\usepackage{amssymb,epic,eepic,epsfig,amsmath,amsbsy,amsmath,bbm}
\usepackage{pb-diagram,lamsarrow,pb-lams,amscd}
\usepackage{graphicx}
\usepackage{texdraw}
\usepackage{url}
\usepackage[bookmarks=true,%
    colorlinks=true,%
    linkcolor=blue,%
    citecolor=blue,%
    filecolor=blue,%
    menucolor=blue,%
    urlcolor=blue,%
    breaklinks=true]{hyperref}

\newtheorem{theorem}{Theorem}[section]
\theoremstyle{definition}
\newtheorem{proposition}[theorem]{Proposition}
\newtheorem{lemma}[theorem]{Lemma}
\newtheorem{definition}[theorem]{Definition}
\newtheorem{remark}[theorem]{Remark}
\newtheorem{corollary}[theorem]{Corollary}
\newtheorem{conjecture}[theorem]{Conjecture}
\newtheorem{question}[theorem]{Question}

\def\printname#1{
    \if\draft y
        \smash{\makebox[0pt]{\hspace{-0.5in}
            \raisebox{8pt}{\tt\tiny #1}}}
    \fi
}

\newcommand{\psdraw}[2]
         {\begin{array}{c} \hspace{-1.3mm}
    \raisebox{-4pt}{\epsfig{figure=draws/#1.eps,width=#2}}
    \hspace{-1.9mm}\end{array}}

\newlength{\standardunitlength}
\catcode`\@=11
\long\def\@makecaption#1#2{%
     \vskip 10pt

\setbox\@tempboxa\hbox{
       \small\sf{\bfcaptionfont #1. }\ignorespaces #2}%
     \ifdim \wd\@tempboxa >\captionwidth {%
         \rightskip=\@captionmargin\leftskip=\@captionmargin
         \unhbox\@tempboxa\par}%
       \else
         \hbox to\hsize{\hfil\box\@tempboxa\hfil}%
     \fi}
\font\bfcaptionfont=cmssbx10 scaled \magstephalf
\newdimen\@captionmargin\@captionmargin=2\parindent
\newdimen\captionwidth\captionwidth=\hsize
\catcode`\@=12

\newcommand{\qbinom}[2]{\text{$\left[\begin{array}{c}#1\\ #2\end{array}
\right]$}}

\def\lbl#1{\label{#1}\printname{#1}}

\def\BN{\mathbb N}
\def\BZ{\mathbb Z}

\def\BQ{\mathbb Q}
\def\BR{\mathbb R}
\def\BC{\mathbb C}

\def\D{\Delta}

\def\a{\alpha}

\def\La{\Lambda}
\def\l{\lambda}

\def\w{\omega}

\def\g{\gamma}

\def\e{\epsilon}
\def\ve{\varepsilon}

\def\d{\delta}
\def\b{\beta}
\def\th{\theta}

\def\mat#1#2#3#4{\left(
\begin{matrix}
 #1 & #2  \\
 #3 & #4
\end{matrix}
\right)}

\def\longto{\longrightarrow}

\def\al{\alpha}

\def\w{\omega}
\def\bk{\mathbf{k}}

\def\k{\kappa}
\def\ev{\mathrm{ev}}

\def\fg{\mathfrak{g}}
\def\SL{\mathrm{SL}}

\def\slim{\mathrm{slim}}
\def\wlim{\mathrm{wlim}}
\def\vol{\mathrm{vol}}
\def\maxdeg{\mathrm{deg}_+}
\def\mindeg{\mathrm{deg}_-}
\def\deg{\mathrm{deg}}
\newcommand{\tH}{\tilde H}

\begin{document}


\title[Asymptotics of the colored Jones function of a knot]{
Asymptotics of the colored Jones function of a knot}

\author{Stavros Garoufalidis}
\address{School of Mathematics \\
         Georgia Institute of Technology \\
         Atlanta, GA 30332-0160, USA \newline
         {\tt \url{http://www.math.gatech.edu/~stavros}}}
\email{stavros@math.gatech.edu}
\author{Thang TQ L\^e}
\address{School of Mathematics \\
         Georgia Institute of Technology \\
         Atlanta, GA 30332-0160, USA \newline
         {\tt \url{http://www.math.gatech.edu/~stavros}}}
\email{letu@math.gatech.edu}

\thanks{The authors were supported in part by National Science Foundation. \\
\newline
1991 {\em Mathematics Classification.} Primary 57N10. Secondary 57M25.
}

\date{August 26, 2011}

\dedicatory{Dedicated to Louis Kauffman on the occasion of his 60th birthday}

\begin{abstract}
To a knot in 3-space, one can associate a sequence of Laurent polynomials,
whose $n$th term is the $n$th colored Jones polynomial. The paper is concerned
with the asymptotic behavior of the value of the $n$th colored Jones polynomial
at $e^{\a/n}$, when $\a$ is a fixed complex number and $n$ tends to infinity.
We analyze this asymptotic behavior to all orders in $1/n$ when $\a$ is a
sufficiently small complex number.
In addition, we give upper bounds for the coefficients and degree of the
$n$th colored Jones polynomial, with applications to upper bounds in the
Generalized Volume Conjecture. Work of Agol-Dunfield-Storm-W.Thurston implies
that our bounds are asymptotically optimal. Moreover, we give results for
the Generalized Volume Conjecture when $\a$ is near $2 \pi i$. Our proofs use
crucially the cyclotomic expansion of the colored Jones function, due to
Habiro.
\end{abstract}

\maketitle

\tableofcontents


\section{Introduction}
\lbl{sec.intro}

\subsection{Asymptotics of the colored Jones function of a knot}
\lbl{sub.asymptotics}

To a knot $K$ in 3-space, one can associate a sequence of Laurent polynomials
$$
J_{K,n}(q) \in \BZ[q^{\pm 1}]
$$
for $n \in \BN=\{1,2,3,\dots\}$. $J_{K,1}(q)=1$ and
$J_{K,2}(q)$ is the famous {\em
Jones polynomial} of $K$ introduced by Jones in \cite{J}, and
$J_{K,n}(q)$ are roughly speaking the Jones polynomials of
$(n-1)$-parallels of the knot. More precisely, $J_{K,n}(q)$ is the
{\em quantum group} invariant of $K$ using the $n$-dimensional
irreducible $\mathfrak{sl}_2(\BC)$ representation, normalized by
$J_{\text{unknot},n}(q)=1$ for all $n$; see \cite{RT,Tu}. The
sequence $\{J_{K,n}(q)\}_n$ is often called the {\em
colored Jones function} of the knot $K$.

The paper is concerned with the asymptotic growth of the colored Jones
function. More precisely, fix a knot $K$ and consider the sequence
of holomorphic functions:

\begin{equation*}
 f_{K,n}:\BC \longto \BC, \quad
f_{K,n}(z):=J_{K,n}(e^{z/n})
\end{equation*}
for $n \in \BN$. In other words, we are evaluating the $n$-th polynomial
$J_{K,n}(q)$ at a complex $n$-th root of $e^{z}$.
We will be concerned with strong and weak
convergence of the sequence $f_{K,n}$, for $n \in \BN$. Let us
explain what we mean by that. Fix an open subset $U$ of $\BC$
containing 0.

\begin{definition}
\lbl{def.strongweak} {\rm(a)} A sequence of holomorphic functions
$f_n: U \longto \BC$ {\em strongly converges in $U$} to a
holomorphic function $f: U \longto \BC$ (and write
$\slim_{n\to\infty} f_n(z)=f(z)$) if $f_n(z)$ converges to $f(z)$
uniformly   on any compact subset
of $U$. \\
{\rm(b)} A sequence of holomorphic functions $f_n: U \longto \BC$
{\em weakly converges} to a holomorphic function $f: U \longto \BC$
(and write $\wlim_{n\to\infty} f_n(z)=f(z)$) if the Taylor series of
$f_n(z)$ at $z=0$ coefficient-wise converges to the Taylor series of
$f(z)$. In other words, for every $k \geq 0$, we have:
$$
\lim_{n\to\infty} \frac{d^k f_n}{dz^k}\Big|_{z=0}
= \frac{d^k f}{dz^k}\Big|_{z=0}.
$$
\end{definition}

It is easy to see that strong convergence of holomorphic functions
implies weak convergence. The converse is not true (see however, Lemma
\ref{lem.complex} below).

The Melvin-Morton-Rozansky (MMR, in short) Conjecture, which was
settled by Bar-Natan and the first author in \cite{B-NG}, compares
the function $f_{K,n}$ of a knot $K$ with the {\em Alexander
polynomial} $\D_K$ of $K$, normalized by $\D_K(t^{-1})= \D_K(t)$ and
$\D_K(1)=1$.

\begin{theorem}{\rm (The MMR conjecture)}
\lbl{thm.MMR}\cite{B-NG}
For every knot $K$ we have
\begin{equation*}
\wlim_{n \to \infty} \,\, f_{K,n}(z) = \frac{1}{\D_K(e^{z})} .
\end{equation*}
\end{theorem}

Our sample result is the following  analytic form of the MMR
Conjecture, which has application in the Generalized Volume
Conjecture.

\begin{theorem}
\lbl{thm.11} {\rm (Proof in Section \ref{sub.complex})}
For every knot $K$ there exists an open neighborhood
$U_K$ of $0 \in \BC$ such that in $U_K$, we have
\begin{equation*}
\slim_{n\to\infty} f_{K,n}(z)=\frac{1}{\D_K(e^{z})}.
\end{equation*}
\end{theorem}

Given Theorem \ref{thm.11} one may ask for a full asymptotic
expansion of  $f_{K,n}(z)$ in terms of powers of $1/n$. In order to
formulate our results, let us introduce the notion of strong and
weak asymptotic expansions.

\begin{definition}
\lbl{def.ass} Fix an open set $U$ of $\BC$, and holomorphic
functions $f_n:U \longto \BC$ and $R_n: U \longto \BC$.
\newline
\rm{(a)} We will say
that the sequence $f_n$ is {\em strongly asymptotic in $U$} to the
series $ \sum_{k=0}^\infty R_k(z) \left(\frac{z}{n}\right)^k $, and
write
\begin{equation}
\lbl{eq.defstrongass}
f_n(z) \sim^s_{n \to \infty}
\sum_{k=0}^\infty
R_k(z) \left(\frac{z}{n}\right)^k
\end{equation}
if for every $N \geq 0$ we have:
\begin{equation}
\lbl{eq.assstrong}
\slim_{n\to\infty}
\left( \frac{n}{z} \right)^{N} \left(
f_n(z)-\sum_{k=0}^{N-1}
R_k(z) \left(\frac{z}{n}\right)^k
\right)=R_{N}(z).
\end{equation}
\rm{(b)}
Likewise, we will say
that the sequence $f_n$ is {\em weakly asymptotic in $U$} to the
series $ \sum_{k=0}^\infty R_k(z) \left(\frac{z}{n}\right)^k $, and
write
\begin{equation}
\lbl{eq.defweakass}
f_n(z) \sim^w_{n \to \infty}
\sum_{k=0}^\infty
R_k(z) \left(\frac{z}{n}\right)^k
\end{equation}
if for every $N \geq 0$ we have:
\begin{equation}
\lbl{eq.assweak}
\wlim_{n\to\infty}
\left( \frac{n}{z} \right)^{N} \left(
f_n(z)-\sum_{k=0}^{N-1}
R_k(z) \left(\frac{z}{n}\right)^k
\right)=R_{N}(z).
\end{equation}
\end{definition}

Usually, sequences of holomorphic functions $f_n(z)$ do not have
asymptotic expansions (or even a limit, as $n \to \infty$). However,
sequences that appear in perturbative expansions of {\em Quantum
Field Theory} are generally expected to have asymptotic expansions.
In fact asymptotic expansions are generally easier to define (via
{\em Feynman diagram} techniques) than the partition functions
$f_{K,n}(z)$ themselves. Even when the partition functions can be
defined, the asymptotic expansions is a numerically useful way to
approximate them.

In \cite{Ro1},
Rozansky discovered that the sequence $f_{K,n}(z)$ has a weak asymptotic
expansion, where the terms are rational functions in the variable $e^z$.
More precisely, Rozansky proved the following result.

\begin{theorem}
\lbl{thm.Zrat}\cite{Ro1}
For every knot $K$ there exists a sequence  $P_{K,k}(q) \in \BQ[q^{\pm 1}]$
of Laurent polynomials with $P_{K,0}(q)=1$ such that
\begin{equation}
\lbl{eq.Zrat}
f_{K,n}(z) \sim_{n \to \infty}^w
\sum_{k=0}^\infty
\frac{P_{K,k}(e^{z})}{\D_K(e^{z})^{2k+1}} \left(\frac{z}{n}\right)^k .
\end{equation}
\end{theorem}

A different proof, valid for all simple Lie groups, was given in \cite{Ga1},
using work of \cite{GK}. Our result is strong version of Theorem
\ref{thm.Zrat}.

\begin{theorem}
\lbl{thm.1} {\rm (Proof in Section \ref{sub.thm1proof})}
For every knot $K$ there exists an open neighborhood
$\tilde U_K$ of $0 \in \BC$ such that in $ \tilde U_K$, we have
\begin{equation}
\lbl{eq.thm1}
f_{K,n}(z) \sim_{n \to \infty}^s \sum_{k=0}^\infty
\frac{P_{K,k}(e^{z})}{\D_K(e^{z})^{2k+1}}
\left(\frac{z}{n}\right)^k.
\end{equation}
\end{theorem}

\subsection{The generalized volume conjecture}
\lbl{sub.volume}

In this section we state some new information about the Volume Conjecture;
the latter connects two very different approaches to knot theory,
namely Topological Quantum Field Theory and Riemannian
(mostly Hyperbolic) Geometry.

\begin{conjecture}
\lbl{conj.vc}\cite{K,MM}
For every hyperbolic knot $K$ in $S^3$ we have:
\begin{equation*}
\lbl{eq.VC}
\lim_{n \to \infty} \frac{\log| f_{K,n}(2 \pi i)|}{n}=
\frac{1}{2 \pi}\, \vol(\rho_{2\pi i}),
\end{equation*}
where $\vol(\rho_{2 \pi i})$ is the {\em hyperbolic volume} of the
 the knot complement $S^3-K$.
\end{conjecture}

In other words, the sequence $f_{K,n}(2 \pi i)$ of complex numbers
grows exponentially with respect to $n$, and the exponential
growth-rate is proportional to the volume of a hyperbolic knot.

One can define the volume function $\vol(\rho)$ of every
representation $\rho: \pi_1(S^3\setminus K)\to \SL_2(2,\BC)$, see
\cite{Du,CCGLS,Th}, and $\vol(\rho_{2\pi i})$ is exactly the value of this
volume function with $\rho_{2\pi i}$ being the discrete faithful
representation of the knot group.

The idea of the Generalized Volume Conjecture (formulated in part by
Gukov in \cite{Gu}) is that we should use other representations of
the knot complement in $\SL(2,\BC)$. For $\a$ nearby $2 \pi i$, in a
small neighborhood of $\rho_{2\pi i}$ there is a unique (up to
conjugation) representation
$$
\rho_{\a}: \pi_1(S^3-K) \longto \SL(2,\BC)
$$
which satisfies
\begin{equation}
\lbl{eq.meridian}
\rho_{\a}(\text{meridian})=
\mat {e^{\a}} {\star} 0 {e^{-\a}} .
\end{equation}
Alas, there is an additional difficulty. Namely, when $\a/(2 \pi i)$ is rational, we should distinguish
two cases: $\a/(2 \pi i)=1 $ or $\a/(2 \pi i)\neq 1$.
The Generalized Volume Conjecture for $\a$ sufficiently close to $2 \pi i$
may now be stated as follows.

\begin{conjecture}
\lbl{conj.gvc}
\lbl{eq.gvc}
If $\a/(2 \pi i) \in (\BR-\BQ) \cup \{1\}$ is sufficiently close to $1$
then
\begin{equation}
\lbl{eq.GVC}
\lim_{n \to \infty} \frac{\log| f_{K,n}(\a)|}{n} =
c_{\a}\, \vol(\rho_{\a}),
\end{equation}
and if $\a/(2 \pi i) \in \BQ- \{1\}$, then
\begin{eqnarray*}
\lbl{eq.GVCa}
\limsup_{n \to \infty} \frac{\log| f_{K,n}(\a)|}{n} &=&
c_{\a}\, \vol(\rho_{\a}),   \\
\lbl{eq.GVCb}
\liminf_{n \to \infty} \frac{\log| f_{K,n}(\a)|}{n} &=& 0,
\end{eqnarray*}
where $c_{\a} \neq 0$ are some nonzero constants.
\end{conjecture}

The distinction of $\a/(2 \pi i)$ being rational or not is a bit
with odds with the notion of {\em hyperbolic Dehn surgery}
developed by Thurston in \cite{Th}. When
$\a/(2 \pi i) $ is a rational number, the hyperbolic Dehn surgery
theorem associates an orbifold filling to the knot complement whose
volume is $\vol(\rho_{\a})$. Orbifolds are mild generalizations of
manifolds. On the other hand, when $\a/(2 \pi i)$ is irrational,
hyperbolic Dehn surgery associates a space which is topologically a
1-point compactification of the knot complement, with volume
$\vol(\rho_{\a})$. In the following, we will refer to the parameter
$\a$ in the Generalized Volume Conjecture as {\em the angle}, making
contact with standard terminology from hyperbolic geometry.

There are two rather independent parts in the Volume Conjecture:
\begin{itemize}
\item[(a)]
To show that the limit exists in \eqref{eq.GVC},
\item[(b)]
To identify the limit with the volume of the corresponding Dehn filling.
\end{itemize}

At the moment, the Generalized Volume Conjecture is known only for the $4_1$
knot and certain values of $\a$; see Murakami, \cite{M1}.

One may further ask what happens to the Generalized Volume
Conjecture when the angle $\a$ is small. For $\a=0$, it is natural
to define $\rho_0$ to be the {\em trivial} representation. Then for
$\a$ small enough, there is a unique (up to conjugation) {\em
abelian} $\SL_2(\BC)$ representation $\rho_{\a}$ that satisfies
\eqref{eq.meridian}.
Abelian representations have $0$ volume (see eg, \cite{CCGLS}).
On the other hand, for small enough
$\a$, we have $\D_K(e^{\a}) \sim \D_K(1)=1$. Thus Theorem \ref{thm.11}
implies that

\begin{theorem}
\lbl{thm.1alt}
For every knot $K$ there exists an open neighborhood $U_K$ of $0 \in \BC$,
such that for $\a \in U_K$, we have:
$$
\lim_{n \to \infty} \frac{\log| f_{K,n}(\a)|}{n}=0=\vol(\rho_{\a}).
$$
\end{theorem}

In other words, Theorem \ref{thm.11} settles the Generalized Volume
Conjecture for small complex angles.

\subsection{The Generalized Volume Conjecture near $2 \pi i$}
\lbl{sub.faithfulresults}

Our next result states that the volume conjecture can only be barely true.

\begin{theorem}
\lbl{thm.near1} {\rm (Proof in Section \ref{sec.cyclotomic})}
For every knot $K$ and every fixed integer $m \neq
0$
$$
\lim_{n\to\infty}  \frac{1}{n} \log|J_{K,n+m}(\exp(2 \pi i/n))| =0.
$$
\end{theorem}

It follows that the {\em double-scaling limit}
$$
\lim_{n,k} \frac{1}{n} \log|J_{K,n}(\exp(2 \pi i/k))|
$$
when $n,k \to \infty$ and $n/k \to 1$ does not exist, or equals to $0$;
with the latter case in contradiction with the Volume Conjecture.
Our next result confirms the strange behavior in the Generalized Volume
Conjecture when $\a/(2 \pi i)$ is rational, not equal to $1$.

\begin{theorem}
\lbl{thm.near2} {\rm (Proof in Section \ref{sec.cyclotomic})}
For every knot $K$ there exist a neighborhood $V_K$
of $1 \in \BC$ such that when  $\a/(2 \pi i) \in V_K$ is rational
and not equal to $1$, then
$$
\liminf_{n\to\infty} \frac{|f_{K,n}(\a)|}{n} =0.
$$
\end{theorem}

\subsection{Upper bounds for the generalized volume conjecture}
\lbl{sub.results}

Our next theorem is an upper bound for the Generalized Volume Conjecture. Let $\Re(\a)$ denote the real part of $\a$.

\begin{theorem}
\lbl{thm.up}
{\rm (Proof in Section \ref{sub.thmup})}
For every knot $K$ with $c+2$ crossings and every $\a \in \BC$, we have
$$
\limsup_{n \to \infty} \frac{\log|f_{K,n}(\a)|}{n} \leq c  \log 4 +
\frac{c+2}{2}|\Re(\a)|.
$$
\end{theorem}

\subsection{Relation with hyperbolic geometry, and asymptotically sharp bounds}
\lbl{sub.hypgeom}

When $\a=2 \pi i$, the upper bound in Theorem \ref{thm.up} is not
optimal, and does not reveal any relationship between the $\limsup$
and hyperbolic geometry. Our next theorem fills this gap.

\begin{theorem}
\lbl{thm.4}
{\rm (Proof in Section \ref{sub.thm4})}
For every knot $K$ with $c+2$ crossings we have
$$
 \limsup_{n \to \infty} \frac{\log|f_{K,n}(2\pi i)|}{n}
\leq  \frac{v_8}{2 \pi} c,
$$
where
$$
v_8=8 \La(\pi/4) \approx 3.6638623767088760602 \dots
$$
is the volume of the {\em regular ideal octahedron}--see \cite{Th}.
\end{theorem}

Using an ideal decomposition of a knot complement by placing one octahedron
per crossing, it follows that for every knot $K$ with $c+2$ crossings, we have
\begin{equation}
\lbl{eq.vol8}
\vol(S^3-K) \leq v_8 c,
\end{equation}
where $\vol(S^3-K)$ is the hyperbolic volume of the knot complement.
On the other hand, if the volume conjecture holds for $\a=2 \pi i$, then
$$
 \lim_{n \to \infty} \frac{\log|f_{K,n}(2\pi i)|}{n}
= \frac{1}{2 \pi} \vol(S^3-K) \leq \frac{v_8}{2 \pi} c.
$$
One may ask whether \eqref{eq.vol8} (and therefore, whether the
bound in Theorem \ref{thm.4}) is optimal. This may be a little surprising,
since it involves all knots (and not just alternating ones)
and their number of crossings, an invariant that carries little known
geometric information.
In conversations with I.Agol and D.Thurston, it was communicated to us
that the upper bound in \eqref{eq.vol8} is indeed optimal. Moreover a class
of knots that achieves (in the limit) the optimal ratio of volume
to number of crossings is obtained by
taking a large chunk of the following {\em weave}, and closing it up to a knot:
$$
\psdraw{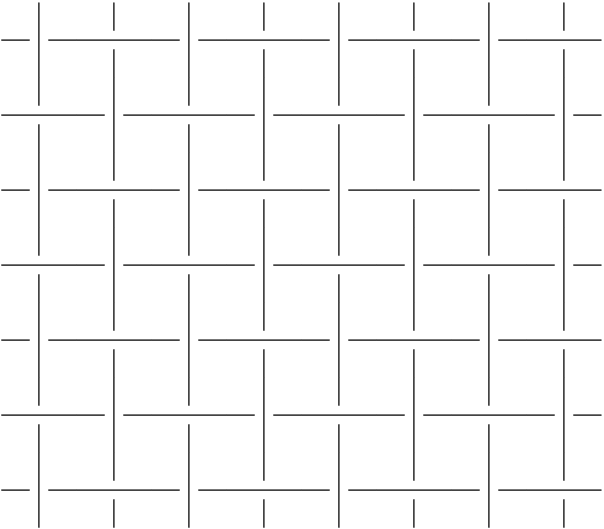}{1.7in}
$$
The complement of the weave has a complete hyperbolic structure associated
with the {\em square tessellation} of the Euclidean plane:
$$
\psdraw{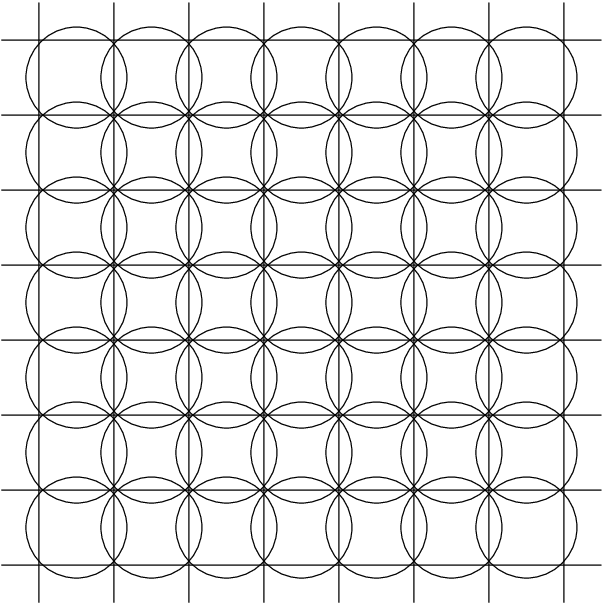}{1.7in}
$$
Optimality follows along similar lines as the Appendix of \cite{La},
using a stronger estimate for the lower bound of the volume of Haken
manifolds, cut along an incompressible surface: If M is a hyperbolic
finite volume 3-manifold containing a properly imbedded orientable,
boundary incompressible, incompressible surface $S$, then
$$
\vol(M) \geq \vol(\mathrm{Guts}(M-\mathrm{int(nbd(}S))),
$$
where $\vol$ stands for volume, and the $\mathrm{Guts}$ terminology
are defined in \cite{Ag}. The proof of this
 stronger statement (of
Agol-Dunfield-Storm-W.Thurston \cite{ADST}) uses, among other
things, work of Perelman.

The reader may compare \eqref{eq.vol8} with the following result of
Agol-Lackenby-D.Thurston \cite{La}:  If $K$ is an alternating knot
with  a planar projection having $t$ twist, then
$$
v_3 (t-1)/2 <  \vol(S^3-K) < 10 v_3 (t-1),
$$
where $v_3=2 \La(\pi/3)\approx 1.01494\dots$ is the volume of the regular ideal
tetrahedron.
Moreover, the class of knots obtained by Dehn filling on the {\em chain link}
has asymptotic ratio of volume by twist number equal to $10 v_3$.
The corresponding tessellation of the Euclidean plane is given by the
{\em star of David}.

So far, we have formulated a Generalized Volume Conjecture for $\a$
near $0$ and $\a$ near $2\pi i$, using representations near the
trivial or near the discrete faithful. How can we
connect these choices for other complex angles $\a$? A natural
answer to this question requires analyzing asymptotics of solutions
of difference equations with a parameter. This is a different
subject that we will not discuss here; instead we will refer the
curious reader to \cite{GG}, and forthcoming work of the first
author. For a further discussion, see also Section
\ref{sec.physics}.

\subsection{The main ideas and organization of the paper}
\lbl{sub.ideas}

In section \ref{sub.complex} we show that weak convergence plus
uniform boundedness implies strong convergence. Thus the strong
convergence of
Theorems \ref{thm.11} and \ref{thm.1} follows from the weak
convergence of Theorems \ref{thm.MMR} and \ref{thm.Zrat},
plus uniform bounds.  Uniform bounds for the colored Jones function
require large cancellations. In order to control these cancellations,
we use the cyclotomic expansion of the colored Jones function of a
knot, which is recalled in Section \ref{sec.fti}. An important point
about this expansion is that its kernel can absorb the exponential
bounds of the coefficients of the cyclotomic functions; see Sections
\ref{sec.estimates} and \ref{sec.allorders}.

Using a state-sum formula for the colored Jones function, we give in
Section \ref{sec.boundsJones} bounds for the degrees and
coefficients of the $n$-th colored Jones polynomial. The result is
also of independent interest. The important point is that the local
weights in the state-sum formula (i.e., the entries of the
$R$-matrix) are Laurent polynomials, given by some ratio of
$q$-factorials. A priori, the bounds of the $n$-colored Jones
function are not good enough to deduce the bounds for the $n$-th
cyclotomic function. However, in Section \ref{sec.boundscyclotomic},
we use a lemma on the growth-rate of the number of partitions of an
integer, in order to deduce the desired bounds for the cyclotomic
function. As a corollary, we can deduce the upper bound of Theorem
\ref{thm.up}.

In the independent Section \ref{sec.hypgeom}, we give a better bound
for the growth-rate of the entries of the $R$-matrix. The important
point is that these entries are ratio of $5$ $q$-factorials, and
each $q$-factorial grows exponentially with rate given by the
Lobachevsky function. The $q$-factorials are arranged in such a way
to deduce that the exponential growth-rate of the entries of the
$R$-matrix is given by the volume of an ideal octahedron. Together
with our state-sum formulas for the $n$-th colored Jones polynomial,
it results in the upper bound of Theorem \ref{thm.4}.

We discuss in Section \ref{sec.cyclotomic} the proof of Theorems
\ref{thm.near1} and \ref{thm.near2}.

In Section \ref{sec.qholonomic} we discuss  bounds on the degrees
and coefficients of $q$-holonomic functions. Earlier work of the authors
implies that the colored Jones and the cyclotomic functions of a knot are
$q$-holonomic.

In Section \ref{sec.physics} we discuss some physics ideas related to the
various expansions of the colored Jones function.

Finally, in  the Appendices we establish the Volume Conjecture for the
Borromean rings using estimates obtained in the proofs of the main results.
At the time when the first draft of this paper was written (2004),  this was
the only hyperbolic link for which the volume conjecture is established.
Since then the volume conjecture has been proved for several other hyperbolic
links, see eg. \cite{Veen}.

Note that Theorems \ref{95} and \ref{thm.L1bound} are not used in the
proofs of our results, and are of independent interest.

The logical dependence of the main theorems is as follows:

$$
\divide\dgARROWLENGTH by2
\begin{diagram}
\node{\text{Thm \ref{thm.1alt}}}
\node{\text{Thm \ref{thm.near1}}}
\node{\text{Thm \ref{thm.near2}}}
\\
\node[2]{\text{Thm \ref{thm.11}}}
\arrow{nw}\arrow{n}\arrow{ne}
\node[2]{\text{Thm \ref{thm.1}}}
\\
\node{\text{Thm \ref{thm.MMR}}}
\arrow{ne}
\node{\text{Thm \ref{thm.boundC3}}}
\arrow{n}
\node{\text{Thm \ref{thm.boundC4}}}
\arrow{ne}
\node{\text{Thm \ref{thm.Zrat}}}
\arrow{n}
\\
\node[2]{\text{Thm \ref{thm.boundC}}}
\arrow{n}\arrow{ne}
\node[2]{\text{Thm \ref{thm.up}}}
\node[2]{\text{Thm \ref{thm.4}}}
\\
\node[3]{\text{Thm \ref{thm.3}}}
\arrow{nw}\arrow{ne}
\node[2]{\text{Prop \ref{prop.R+limit}}}
\arrow{ne}
\node[2]{\text{Prop \ref{prop.R}}}
\arrow{nw}
\end{diagram}
$$

\subsection{Acknowledgment}

The authors wish to thank I. Agol, D. Boyd, N. Dunfield, D.
Thurston and D. Zeilberger for many enlightening conversations, and
the anonymous referee for his careful reading of our manuscript.

\section{Weak versus strong convergence}
\lbl{sec.versus}

\subsection{A lemma from complex analysis}
\lbl{sub.complex}

To prove Theorem \ref{thm.11}, we need to improve the weak
convergence of Theorem \ref{thm.MMR} to the strong convergence.
This uses the next lemma on normal
families that is sometimes referred to by the name of {\em Vitali} or
{\em Montel's} theorem. For a reference, see \cite{Hi,Sch}. The lemma exhibits
the power of holomorphy, coupled with uniform boundedness.


\begin{lemma}
\lbl{lem.complex}
If
$$
f_n: \{z \in \BC \, : \, |z| < r \} \to \{z \in \BC \, : \, |z| \leq M \}
$$
is a uniformly bounded sequence of holomorphic functions such
that for every $m \geq 0$, we have:
$$
\lim_{n \to \infty} f^{(m)}_n(0) =a_m.
$$
Then,
\begin{itemize}
\item
The limit $f(z)=\lim_n f_n(z)$ exists pointwise for all $z$ with $|z| < r$ .
\item
$f$ is holomorphic,
\item
The convergence is uniform on compact subsets, and
\item
For every $m$, $f^{(m)}(0)=a_m$.
\end{itemize}
In other words, weak convergence and uniform boundedness imply
strong convergence.
\end{lemma}

\begin{proof}
$\{f_n\}_n$ is uniformly bounded, so it is a normal family, and contains
a convergent subsequence $f_j\to f$. Convergence is uniform on compact sets,
and $f$ is holomorphic, and for every $m \geq 0$,
$\lim_j f_j^{(m)}(0)=f^{(m)}(0)=a_m$.

If $\{f_n\}_n$ is not convergent (uniformly on compact sets),
since it is a normal family, then
there exist two subsequences that converge to $f$ and $g$
respectively, with $f \neq g$. Applying the above discussion, it
follows that $f$ and $g$ are holomorphic functions with equal
derivatives of all orders at $0$. Thus, $f=g$, giving a
contradiction.
\end{proof}

 Theorem \ref{thm.11} follows from Lemma \ref{lem.complex} and
the following result, whose proof will be given in Section
\ref{sec.estimates}.

\begin{theorem}
\lbl{thm.boundC3} {\rm (Proof in Section \ref{sub.thm.boundC2})}
For every knot $K$ there exists an open neighborhood $U_K$ of $0 \in
\BC$ and a positive number $M$ such that for $\a \in U_K$, and all
$n \geq 1$, we have:
$$
|f_{K,n}(\a)| < M.
$$
\end{theorem}

\subsection{The main difficulty for uniform bounds}
\lbl{sub.keydiff}

Before we proceed with the proof of Theorem \ref{thm.boundC3}, let us
point out the main difficulty. As we will see later, $J_{K,n}(q)$
is a Laurent polynomial in $q$ whose span (i.e., the exponents
of its monomials)
are $O(n^2)$ and whose coefficients are $e^{O(n)}$. In addition, due to
our normalization, $J_{K,n}(1)=1$. In other words, the $O(n^2)$ many
exponentially growing
coefficients of $J_{K,n}(q)$ add up to $1$. When we evaluate
$J_{K,n}(e^{\a/n})$, we want to bound the result independent of $n$. This
will happen only if major cancellations occur. How can we control these
cancellations? The answer to this is a key cyclotomic expansion of the
colored Jones function, which we review next.

\section{Two expansions of the colored Jones polynomial}
\lbl{sec.fti}

\subsection{The loop expansion}
\lbl{sub.loop}

With $q=e^h$, one has

$$
J_{K,n}(e^h) = \sum_{i=0}^\infty a_{K,i}(n)\,
h^i \in \BQ[[h]].
$$

It turns out that $a_{K,i}(n)$ is a polynomial in $n$ with degree less
than or equal to $i$, see \cite{B-NG}. Hence there are rational numbers
$a_{K,i,j}$, depending on the knot $K$, such that

\begin{eqnarray}
\lbl{eq.pert} J_{K,n}(e^h) &=& \sum_{0 \leq j \leq i} a_{K,i,j}\, n^j
h^i  = \sum_{0 \leq i, 0 \leq j \leq i} a_{K,i,j} (nh)^j h^{i-j}\notag \\
\notag &=& \sum_{0 \leq j, k} a_{K,j+k,j} (nh)^j h^k.
\end{eqnarray}
If we define
$$
R_{K,k}(x)=\sum_{0 \leq j} a_{K,j+k,j} x^j \in \BQ[[x]],
$$
then we have the following {\em loop expansion}
\begin{eqnarray}
\lbl{eq.R}
J_{K,n}(e^h) &=& \sum_{k=0}^\infty R_{K,k}(nh) h^k \in \BQ[[h]]
\end{eqnarray}

It turns out that $R_{K,k}(x) \in \BQ(e^x)$ are rational functions
for all $k$. In fact,
the MMR Conjecture states that
$$
R_{K,0}(x)=\frac{1}{\D_K(e^x)} \in \BQ[[x]].
$$
More generally, Rozansky \cite{Ro1} proves  there are Laurent
polynomials $P_{K,k}(t) \in \BQ[t^{\pm 1}]$ such that in $\BQ[[x]]$,
$$
R_{K,k}(x)=\frac{P_{K,k}(e^x)}{\D_K(e^x)^{2k+1}}.
$$

\begin{remark}
\lbl{rem.aij}
For every $i,j$, the function $K \to a_{K,i,j}$ is a finite type invariant
of degree $i$.
Although the polynomials $P_{K,k}(t)$ are not finite type invariants
(with respect to the usual crossing change of knots), they are
finite type invariants with respect to a loop move described in \cite{GR}.
We will not use these facts in our paper.
\end{remark}

\subsection{The cyclotomic expansion}
\lbl{sub.cyclotomic}

Habiro  found another interesting expansion of the colored Jones function, known as the
cyclotomic expansion. Although the cyclotomic expansion has important
arithmetic consequences, we discuss only its algebraic properties here. Let
us define:

\begin{equation}
\lbl{eq.Ckernel}
C_{n,k}(q) = \prod_{j=1}^{k} (q^n+q^{-n}-q^j-q^{-j}),
\quad \text{with} \quad C_{n,0}(q)
 :=1.
\end{equation}

Habiro  showed that there exist unique Laurent polynomials $H_{K,k}(q) \in \BZ[q^{\pm 1}], k=0,1,\dots $ such that

\begin{equation}
\lbl{cyc}
J_{K,n}(q) = \sum_{k=0}^{n-1} C_{n,k}(q)\, H_{K,k}(q).
\end{equation}

For details, see \cite[Section 6]{H1}. Note that our $H_{K,n}(q)$ is  $J_K(P''_n)$ in  Habiro's notation.
We will call the expansion
\eqref{cyc} the cyclotomic expansion. Since $C_{n,k}(q)=0$ if $k
\ge n$, the summation in \eqref{cyc} can be assumed from $0$
to $\infty$.


It is possible to solve for $H_{K,n}$ from Equation \eqref{cyc}. Explicitly, from \cite[Lemma 6.1]{H1} one has
\begin{equation}
  H_{K,n}(q)=\frac{1}{\{2n+2\}!}\sum_{k=1}^{n+1} (-1)^{n+1-k}
  \{2k\}\{k\}
\qbinom{2n+2}{n+1-k} J_{K,k}(q) \label{33}
\end{equation}

where we use the following definition

$$ \{n\} := q^{n/2} -q^{-n/2} \quad \text{and} \quad \{n\}! :=
\prod_{i=1}^n \{i\},$$

\begin{equation*}
 \qquad \{a\}_b :=
\frac{\{a\}!}{\{a-b\}!}=\prod_{j=a-b+1}^a\{ j \}, \qquad
\qbinom{a}{b} := \frac{\{ a \}!}{\{ b \}! \{ a-b \}!}
\end{equation*}

\subsection{Comparing the cyclotomic and the loop expansion}
\lbl{sub.comparing}

In the loop expansion, as well as in the
cyclotomic expansion, one should treat $q^n$ and $q$ (where $n$ is
the color) as two independent variables. Consider two
independent variables $z$ (standing for $\a$) and $y$ (standing for
$\a/n$). Let us define the following biholomorphic functions

\begin{eqnarray*}
c_k(z,y) &=&
\prod_{j=1}^{k} (e^z+e^{-z} -e^{jy}-e^{-jy}), \\
h_{K,k}(z,y) &=& c_k(z,y) H_{K,k}(e^{y}).
\end{eqnarray*}

The cyclotomic expansion says that for every $n$ we have:

\begin{equation}
\lbl{11}
f_{K,n}(\a) = \sum_{k=0}^{\infty} h_{K,k}(\a,\a/n) \in \BQ[[\a]].
\end{equation}

The loop expansion is a Taylor expansion in $\a/n$, so we will
consider the Taylor expansion in $y$ (around 0) of $h_{K,k}(z,y)$:

$$
h_{K,k}(z,y) = \sum_{p=0}^\infty d_{k,p}(z) \, y^p,
$$
where $d_{k,p}(z)$ (which depends on $K$) is holomorphic for $z \in \BC$.

Comparing the loop and the cyclotomic expansion (Equations \eqref{eq.R}
and \eqref{11}), we obtain that:

\begin{lemma}
\lbl{lem.compare3}
For every knot $K$ and every $p \in \BN$ we have
\begin{equation} R_{K,p}(x) = \sum_{k=0}^\infty
d_{k,p}(x) \in \BQ[[x]].\label{12}
\end{equation}
as formal power series in $x$.
\end{lemma}

\section{A reduction of Theorem \ref{thm.boundC3}
to estimates of the cyclotomic function}
\lbl{sec.estimates}

\subsection{Uniform bounds of the colored Jones function}
\lbl{sub.reductionA}

In this section we will deduce Theorem \ref{thm.boundC3} from estimates
of the degree and the coefficients of the cyclotomic expansion of
the knot. These estimates will be established in Section
\ref{sec.boundsJones}.
By definition $f_{K,n}(\a) = J_{K,n}(e^{\a/n})$, hence equation
\eqref{cyc} gives that
\begin{eqnarray}
\lbl{eq.J2Ca} f_{K,n}(\a) &=& \sum_{k=0}^{n-1} C_{n,k}(e^{\a/n})
H_{K,k}(e^{\a/n}).
\end{eqnarray}

To have upper bounds for $|f_{K,n}(\a)|$ we will need bounds for
$H_{K,k}(e^{\a/n})$ and the ``kernel'' $C_{n,k}(e^{\a/n})$ (the
kernel does not depend on the knot $K$).

\begin{definition}
\lbl{def.l1norm}
For a Laurent polynomial $f(q)=\sum_k a_k q^k$, we define its
$l^1$-norm by
$$
||f||_1=\sum_k |a_k|.
$$
\end{definition}
The proof of the following Theorem, which gives bounds for the
degrees and the $l^1$-norm of $H_{K,n}$, will be given in Section
\ref{sec.boundscyclotomic}.

\begin{theorem}
\lbl{thm.boundC}
{\rm (Proof in Section \ref{sec.boundscyclotomic})}
For every knot $K$, there are positive constants
$A_0, A_1$ (depending on $K$) such that for all $n \in \BN$ we have
\begin{equation}
H_{K,n}(q) = \sum_{j=-A_0n^2}^{A_0n^2} b_{j,n} q^j.
\tag{a}\label{eq1}
\end{equation}
and
\begin{equation}
\lbl{eq.LC} ||H_{K,n}||_1 \leq A_1^n . \tag{b}
\end{equation}
\end{theorem}

The next lemma follows from Theorem \ref{thm.boundC} and an elementary estimate.

\begin{lemma}
\lbl{lem.estimate}
Suppose $|\a| < 1$.
\begin{itemize}
\item[(a)]   For every knot $K$ there is a constant
$A_2$ such that for every $ 0\le k\le n$, we have
$$
|H_{K,k}(e^{\a/n}) | \le (A_2)^k.
$$
\item[(b)]
There is a constant $A_3 >0$ such that every $ 0\le k \leq n$ we
have:
$$
|C_{n,k}(e^{\a/n})| \leq (A_3)^k |\a|^k.
$$
\end{itemize}
\end{lemma}

\begin{proof}
(a) By Theorem \ref{thm.boundC}(a),
$$
H_{K,k}(e^{\a/n}) = \sum_{j=-A_0k^2}^{A_0k^2} b_{j,k} e^{j\a/n}.
$$
From the bounds for $j$ and $k\le n$  one has that $|j/n| \le A_0 k$,
hence $ |e^{j\a/n}| \le \exp (A_0\,  k\, |\Re(\a)|) \le \exp(k \, A_0) $. From
the above equation one has
$$
\left | H_{K,k}(e^{\a/n})\right | \le ||H_{K,k}||_1 \, (\exp A_0)^k.
$$
Using Theorem \ref{thm.boundC}, it is enough to take $A_2 = A_1
\exp(A_0)$.

(b) By definition,
$$
C_{n,k}(e^{\a/n})=\prod_{j=1}^{k} (e^{\a}+e^{-\a}-e^{j\a/n}
-e^{-j\a/n}) = \prod_{j=1}^{k} \left(g(\a) -g(j\a/n)\right),
$$
where $g(z)= e^z + e^{-z}$. One has $g'(z) = e^z -e^{-z}$, hence for
$z$ on the interval connecting $\a$ and $j\a/n$, with $0\le j\le n$,
one has $|g'(z)| \le 2 \exp(|\a|) \le  2e$. By the mean value theorem,
we have, for $0 \le j \le k \le n$,
$$
|g(\a) -g(j\a/n)| \le 2 e |\a -j\a/n| \le  2 e |\a|.
$$
It follows that
$$
|C_{n,k}(e^{\a/n})| \le (2 e)^k |\a|^k.$$ It is enough to take
$A_3=2e$.
\end{proof}

\subsection{Theorem \ref{thm.boundC} implies Theorem \ref{thm.boundC3}}
\lbl{sub.thm.boundC2}

It follows from Lemma \ref{lem.estimate}  that for  $0 \leq k \leq n$
and $|\a|< 1$, we have:
\begin{eqnarray*}
|C_{n,k}(e^{\a/n}) H_{K,k}(e^{\a/n})| & \leq & |\a A_2 A_3|^k.
\end{eqnarray*}

Let us choose $U_K$ to be the disk centered  at the 0, with radius
$1/(2A_2A_3+1)$, then $ |\a A_2 A_3| < 1/2$ for $\a \in U_K$.
Equation \eqref{eq.J2Ca} and the above estimate imply that for all
$n$ and all $\a \in U_K$, we have:
$$
|f_{K,n}(\a)| \le  \sum_{k=0}^{n-1}
\left |C_{n,k}(e^{\a/n}) H_{K,k}(e^{\a/n})\right| \\
\le  \sum_{k=0}^{n-1} (1/2)^k < 2.
$$
which concludes the proof of Theorem \ref{thm.boundC3}, assuming
Theorem \ref{thm.boundC}.
\qed

\section{A reduction of Theorem \ref{thm.1} to estimates of the cyclotomic
function} \lbl{sec.allorders}

\subsection{Some estimates}

The following is a higher order version of Lemma \ref{lem.estimate}.
The proof is similar.

\begin{lemma}
\lbl{lem.estimate2}
Suppose $|\a| < 1$.
\begin{itemize}
\item[(a)]  For $1\le k \le n$, $0\le l$,  and $y$ on
the interval from 0 to $\a/n$ we have
$$
\left |\frac{\partial^l}{\partial y^l} c_k(\a,y)
\right| < (A_3)^k  |\a|^{k-l}\, k^{2l}.
$$
\item[(b)]
For any $y \in \BC, |y| < 1/n$ and $1\le k\le n$ we have
$$
\left |\frac{\partial^l}{\partial y^l}H_{K,k}(e^y)\right
| <(A_2)^k\, (A_0)^l \, k^{2l}.
$$
\end{itemize}
\end{lemma}

\begin{proof}
(a) We have $c_k(\a,y) = \prod_{j=1}^k g_j$, where
$$
g_j = e^\a + e^{-\a} -e^{jy} - e^{-jy}.
$$
\newcommand{\bl}{{\mathbf l}}
By the Leibniz rule, the $l$-th derivative (with
respect to $y$) of $c_k$ is the sum
\begin{equation}
\lbl{form}
\frac{\partial^l}{\partial y^l} c_k(\a,y) = \sum_{|\bl|=l} \binom {l}{l_1,\dots,l_k}\, t(\bl),
\quad \text{where} \quad
t(\bl)=\prod_{j=1}^k g_j^{(l_j)}.
\end{equation}

Here $\bl=(l_1,\dots,l_k), |\bl|:= \sum_{j=1}^k l_j$, $l_j \ge 0$.
 We will estimate each term $t(\bl)$. Fix $\bl= (l_1,\dots,l_k)$. We
consider two cases, $l_j=0$ and $l_j>0 $.

Suppose $l_j=0$. Then $g_j^{(l_j)} = g_j = e^\a + e^{-\a} - e^{jy}
- e^{-jy}$. Since $j \le k \le n$, the interval connecting $\a$ and
$jy$ lies totally in the disk of radius $|\a|$ (remember that $|y|
\le |\a|/n$). As in the proof of Lemma \ref{lem.estimate}, we have
\begin{equation}
\lbl{1}
|g_j|=|(e^{\a} + e^{-\a}) - (e^{jy} + e^{-jy})| \le (2e)|\a|.
\end{equation}

Now suppose $l_j >0$. Then
$$
g_j^{(l_j)} =  - e^{jy}j^{l_j} - e^{-jy}(-j)^{l_j}.$$

It is clear $|e^{\pm \a}| < e$. Since $|j|\le |k|$ and $|jy| < |\a|$,
we have $|e^{\pm jy}(\pm j)^{l_j}| < e k^{l_j}$. Hence

\begin{equation}
\left |g_j^{(l_j)}\right| <  (2e)\, k^{l_j}.
\label{2}
\end{equation}
Taking the product over $j$, using \eqref{1}, \eqref{2} and $\sum l_j =l$,
we get

$$
\begin{aligned} |t(\bl)| & < (2e)^k\,  k^l\,  |\a|^{\#\{l_j =0 \}} \\
        & < (2e)^k\,  k^l\,  |\a|^{k-l} \quad
\text{because} \quad |\al| <1 \quad \text{and} \quad \#\{l_j =0 \} \ge k-l
\end{aligned}
$$
Since $\sum_{|\bl|=l} \binom {l}{l_1,\dots,l_k} =k^l$, from \eqref{form}  and the above estimate for $t(\bl)$, we get the result with
$A_3=2e$.

(b) By Theorem \ref{thm.boundC}(a),

$$
\frac{\partial^l}{\partial y^l}H_{K,k}(e^{y}) =
\sum_{j=-A_0k^2}^{A_0k^2} b_{j,k} e^{jy} j^l.
$$
From the bounds for $j$ and $k\le n$  one has $ |e^{jy}| \le \exp
(A_0 k )$ and $|j^l | \leq (A_0k^2)^l$. From the above equation one has

$$
\left| \frac{\partial^l}{\partial y^l} H_{K,k}(e^{y})\right| \le ||H_{K,k}||_1
\, (\exp A_0)^k (A_0k^2)^l.
$$
Using Theorem \ref{thm.boundC}, it is enough to take $A_2 = A_1
\exp(A_0)$.
\end{proof}

\begin{corollary}
\lbl{14}
For every knot $K$ there are positive constants $A_4, A_5$ such
that
\begin{itemize}
\item[(a)] for $0\le k$, $0\le N$,  and  $|\a| <1$ and $y$  on the
interval from $0$ to $\a/k$, we have
$$ \left |   \frac{\partial^N}{\partial y^N}h_{K,k}(\a,y) \right |
< |\a A_4|^{k-N} (A_5 k^2)^N.
$$
\item[(b)] for $0\le k$, $0\le N$, and  $|\a| <1$  and every positive
integer $n$, we have
$$ \left | h_{K,k}(\a, \a/n) - \sum_{p=0}^{N-1} d_{k,p}(\a) \,
 (\a/n)^p \right|  < \frac{1}{N!}
\left( \frac{\a}{n} \right)^N \, |\a A_4|^{k-N} (A_5 k^2)^N.
$$
\end{itemize}
\end{corollary}

\begin{proof} (a) The $N$-th derivative of $h_{K,k}(\a,y)$, which is
the product of $c_k(\a,y)$ and $H_{K,k}(e^y)$, is the sum of $2^N$
terms, each  of the form

$$
\frac{\partial^l}{\partial y^l} c_k(\a,y) \, \frac{\partial^{N-l}}{\partial
y^{N-l}}H_{K,k}(e^y).
$$

Using  Lemma \ref{lem.estimate2}, the absolute value of the above
term is bounded by $|\a|^{k-l} (A_2A_3)^k (A_0)^{N-l} k^{2N}$, which,
in turn, is less than  $|\a|^{k-N} (A_2A_3)^k (A_0)^{N} k^{2N}$. Hence,
multiplied by $2^N$ we get

$$ \left| \frac{\partial^l}{\partial y^l} h_{K,k}(\a,y) \right| < 2^N
\times |\a|^{k-N} (A_2A_3)^k (A_0)^{N} k^{2N} = (\a A_2A_3)^{k-N}
(2A_0A_2A_3 k^2)^N.
$$
It is enough to take $A_4=A_2A_3$ and $A_5= 2A_0A_2A_3$.

(b) By Taylor's Theorem,
$$
\left | h_{K,k}(\a, \a/n) - \sum_{p=0}^{N-1} d_{k,p}(\a) \,
 (\a/n)^p \right|  < \frac{1}{N!} \left( \frac{\a}{n} \right)^N
 \,\max \left|\frac{\partial^N}{\partial y^N}h_{K,k}(\a,y)\right
 |,
$$
where $\max$ is taken when $y$ is on the interval connecting 0 and
$\a/n$.
Using the  estimate of part (a), we get the result.
\end{proof}

\subsection{Theorem  \ref{thm.boundC} implies Theorem \ref{thm.1}}
\lbl{sub.thm1proof}

To simplify notation, let us define, for a knot $K$,
\begin{equation}
\lbl{eq.simplifynot}
f_{K,n}^{[N]}(z):=J_{K,n}(e^{z/n})-\sum_{k=0}^{N-1}
\frac{P_{K,k}(e^{z})}{\D_K(e^{z})^{2k+1}}
\left(\frac{z}{n}\right)^k.
\end{equation}

Theorem \ref{thm.1} follows from Theorem \ref{thm.Zrat},  Lemma
\ref{lem.complex} and the following uniform bound.

\begin{theorem}
\lbl{thm.boundC4}
For every knot $K$ there exists an open neighborhood $\tilde U_K$ of
$0 \in \BC$ such that for every $N \geq 0$ there exists a positive
number $M_N$ such that for $\a \in \tilde U_K$, and all $n \geq 0$, we
have:
$$
\left|\left( \frac{n}{\a} \right)^N
f_{K,n}^{(N)}(\a) \right| < M_N.
$$
\end{theorem}

\begin{proof}(of Theorem \ref{thm.boundC4}, assuming Theorem
\ref{thm.boundC})
We have the following identities, where the second follows from
\eqref{11} and \eqref{12}:

$$\begin{aligned}
f_{K,n}^{[N]}(\a) & = f_{K,n}(\a) -
\sum_{p=0}^{N-1}
R_{K,p}(\a) \, \left (\frac{\a}{n} \right)^p \\
&= \sum_{k=0}^\infty h_{K,k}(\a,\frac{\a}{n}) - \sum_{p=0}^{N-1}
\left( \sum_{k=0}^\infty d_{k,p}(\a) \right) \left (\frac{\a}{n}
\right)^p
\\ &= \sum_{k=0}^\infty \left [h_{K,k}(\a,\frac{\a}{n}) - \sum_{p=0}^{N-1}
d_{k,p}(\a) \left (\frac{\a}{n} \right)^p \right]
\end{aligned}$$

Using the estimate in Corollary \ref{14}, we see that

$$ \left |\left (\frac{n}{\a}\right)^N f_{K,n}^{[N]}(\a)\right| < \frac{1}{N!}
\sum_{k=0}^\infty \, |\a A_4|^{k-N}\,  (A_5 k^2)^N.
$$
 If
$|\a A_4 | < 1$, the series of the right hand side is absolutely
convergent. It is enough to take $\tilde U_K$ to be the disk centered
at 0 with radius $1/(2A_4+1)$. This proves Theorem
\ref{thm.boundC4}, assuming Theorem \ref{thm.boundC}.
\end{proof}

\section{Bounds for the degree and coefficients of the colored
Jones function}
\lbl{sec.boundsJones}

\subsection{Bounds for the degree}
\lbl{sub.bdegree}

In this section we give a bound for the coefficients of the colored
Jones polynomial, and deduce Theorem \ref{thm.up}. This and the
next section are logically independent from the previous Sections
\ref{sec.estimates} and \ref{sec.allorders}.

For a Laurent polynomial $f(q)=\sum_{k=m}^M a_k q^k$, with $a_m a_M
\neq 0$, let us define $\maxdeg(f)=M$ and $\mindeg(f)=m$. In
\cite{Le} the second author showed that there are quadratic bounds for the
 degrees  of the colored Jones polynomial.

Suppose the knot $K$ has a planar projection
with $c+2$ crossings. Let $\w$ be the writhe number, i.e. the number
of positive crossing minus the number of negative ones.
Then by  \cite[Proposition 2.1]{Le}, taking
into account the change of variable, the framing, and the
normalization, one has the following bounds for the degrees of
$J_{K,n}(q)$.

\begin{proposition}
\lbl{21}
With the above notations, there are constants $s_\pm$ such that
$$
\begin{aligned} \maxdeg(J_{K,n})  & \le \frac{(c+2)(n-1)^2 +
2(n-1)(s_+-1) -\w(n^2-1)}{4}  \\
 \mindeg(J_{K,n}) &\ge -\frac{(c+2)(n-1)^2 + 2(n-1)(s_--1) +\w(n^2-1)}{4}.
\end{aligned}
$$
\end{proposition}
The constants $s_\pm$ have transparent geometric meaning, but we don't need
their exact values here.

Another proof of the quadratic bounds, though less as explicit,  for the
degrees of the colored Jones polynomial using the theory of $q$-holomorphic
functions is given in Section \ref{sub.qholo}.

\subsection{Bounds for the coefficients}
\lbl{sub.coeffbounds}

For a Laurent polynomial $f\in \BZ[q^{\pm 1/4}]$, we define $||f||_1$ as in
Definition \ref{def.l1norm}, i.e. $||f||_1$ is the sum of the absolute
values of its coefficients. Observe that

\begin{equation}
 ||f+g||_1 \leq ||f||_1+||g||_1, \qquad ||fg||_1 \leq
||f||_1 \, ||g||_1. \label{eq.ineq}
\end{equation}

Since $||\{j\}||_1 =2$, we have, for $k \leq n$
\begin{equation}
\left| \left|\{a\}_k\right|\right|_1  = \left|\! \left| \prod_{j=a-k+1}^a\{ j \}\right|\!\right|_1  \leq 2^k \leq 2^n
\label{32}
\end{equation}

It is known that the quantum binomial $\qbinom{m}{k}$ is a Laurent
polynomial in $q^{1/2}$ with {\em positive} integer coefficients,
hence its $l^1$-norm is obtained by putting $q^{1/2}=1$, which is
the classical  binomial $\binom{m}{k}$. One has, if $m \leq
n$,
\begin{equation} \left|\! \left|\qbinom{m}{k}\right|\!\right|_1=\binom{m}{k} \leq 2^m \leq
2^n \label{31}
\end{equation}

\begin{theorem}
\lbl{thm.3}
For every knot $K$ of $c+2$ crossings and
every $n$ we have:

\begin{equation}
\lbl{eq.l1}
||J_{K,n}||_1 \leq n^c 4^{c n}.
\end{equation}
\end{theorem}

\begin{proof}
The proof of the Theorem is easy using the state sum
definition of the colored Jones polynomial: The colored Jones
polynomial is the sum, over all states, of the weights of the
states. There are $n^c$ states, the weight of each is the product of
several $q$-factorials and $q$-binomial coefficients for which an
upper bound can be easily found. Let us now go to the details of the proof.

The knot $K$ is the closure of a $(1,1)$-tangle $T$ (or long knot),
with orientation given by the direction from the bottom boundary
point to the top boundary point. The crossing points of the diagram
of $T$ (on the standard 2-plane) break $T$ into $2c+5$ arcs, two of
which are {\em boundary} (i.e. each contains a boundary point of
$T$). The two crossings adjacent to the boundary arcs are called
{\em boundary crossings}.

To get from $c+2$ to $c$ in the estimate, we will choose
 the $(1,1)$-tangle $T$  such that (1) when going along $T$, starting at
 the bottom boundary point, we must pass the very first  crossing
(resp. very last crossing) by an
 overpass (respectively, an underpass) and (2) the two strands
 at each crossing are pointing upwards, as in the following figure:
\begin{equation}
\lbl{cross}
\psdraw{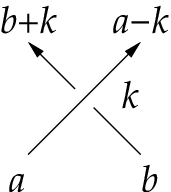}{0.5in} \qquad \psdraw{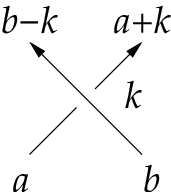}{0.5in}
\end{equation}

Here is how to get such a $(1,1)$-tangle $T$. Consider a diagram of
$K$ on a 2-sphere $S^2$. The $c+2$ crossings break the knot diagram
into $2c+4$ arcs. At each crossing we have an overpass and an
underpass. When we go along the knot starting at some point,
following the direction of the orientation, we pass through all
these underpasses and overpasses. Hence there must be an arc which
starts at an underpass and ends at an overpass, assuming there is at
leat one crossing. Remove from $S^2$ a small disk which is a small
neighborhood of a point inside this arc. What is left is a long knot
diagram on a disk, which can also be considered as a $(1,1)$-tangle
diagram in the strip $\BR\times [0,1]$ in the standard 2-plane which
satisfies requirement (1). Using the isotopy of the form
$$
\psdraw{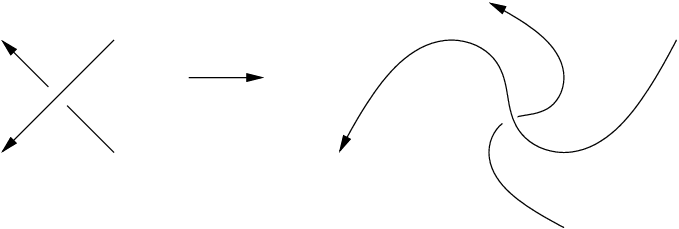}{2in}
$$
which moves crossings (positive or negative) into standard upright
position, we get the desired $(1,1)$-tangle.

A state $\bk$ is an assignment of numbers, called the colors, to the
crossings of the diagram of $T$, where each color is in
$\{0,\dots,n-1\}$. For a fixed state we will color the $2c+5$ arcs
as follow. First color the bottom boundary arc by 0. Going along the
diagram of $T$ from the bottom boundary point, if we are on an arc
of color $a$ and pass a crossing, the next arc will have color $a+k$
or $a-k$, according as the pass is an underpass or an overpass, see
 \eqref{cross}. Here $k$ is the color of the crossing.

We will only consider states such that  the colors of arcs are
between 0 and $n-1$ and the color of the top boundary arc is 0. The
under/overpass configuration at the two boundary crossings ensures
that the two boundary crossings have color 0, otherwise the arcs
next to the two boundary arcs would have negative colors. It follows
that the number of states is at most $n^c$.

The weights of the positive crossing  (on the left) and negative
crossing (on the right in \eqref{cross}) are

\begin{eqnarray}
\lbl{eq.R+-} R_+(n;a,b,k) &=&  (\text{unit})\,
\qbinom{b+k}{k} \{n-1+k-a\}_k, \\
R_- (n;a,b,k) &=& (\text{unit})\,\qbinom{a+k}{k} \{n-1+k-b\}_k,
\end{eqnarray}
where $(\text{unit})$ stands for $\pm$ a power of $q^{\pm 1/4}$,
which does not affect the $l^1$ norm. Note that both $a+k$ and $b+k$
in the above formulas are between $0$ and $n-1$.

The weight of a maximum/minimum point is a also $\pm$  a power of
$q^{\pm 1/4}$, whose exact formula is not important for us.  Let
$F(n,\bk)$ denote the product of weights of all the crossings and
all the extreme points. Then

\begin{equation}
\lbl{eq.statesum}
J_{K,n}(q)=\sum_{\bk} F(n,\bk).
\end{equation}

 Using the estimates \eqref{31} and \eqref{32}, we see
that $||R_{\pm}(n;a,b,k)||_1 \leq 4^n$. Since the weight of the two
boundary crossing is just a unit, the $l^1$ norm of $F(n,\bk)$ is
less than $4^{cn}$. From \eqref{eq.statesum} and the fact that there
are $n^c$ states, we get  $ ||J_{K,n}||_1 \le n^c 4^{cn}$.
\end{proof}

Since there is a constant $b$ such that $n^c \le b^n$, we have the
following.

\begin{theorem}
\lbl{95}
For every knot $K$, there is a constant $A_6$ such
that for every positive integer $n$,
$$ ||J_{K,n}||_1 \le (A_6)^n.$$
\end{theorem}

\subsection{Proof of Theorem \ref{thm.up}}
\lbl{sub.thmup}


Fix a knot with $c+2$ crossings. The bounds for the degrees of
$J_{K,n}$ (see Proposition \ref{21}) allow us to write
$$
J_{K,n}(q)=\sum_{j} a_{n,j} q^j
$$
where $ |j|\leq  n^2 (c+2+|w|)/4 + O(n)$. For such $j$, we have:
$$
|e^{j\a/n}|=e^{\Re(j\a)/n} \leq e^{ (c+2+|w|)n/4 + O(1))|\Re(\a)|}.
$$
Using Theorem \ref{thm.3} we get
$$
|J_{K,n}(e^{\a/n})| \leq  n^c 4^{cn} e^{ (c+2+|w|)n/4 +
O(1))|\Re(\a)|}.
$$
Thus,
$$
\frac{1}{n} \log|f_{K,n}(\a)| \leq c\log4 + \frac{c+2+|w|}{4}
|\Re(\a)| + O\left(\frac{\log n}{n} \right).
$$
The result follows from the observation that $|\w| \le c+2$, since
$c+2$ is the total number of crossings.
\qed

\section{Proof of Theorem \ref{thm.boundC}}
\lbl{sec.boundscyclotomic}

The goal of this Section is to prove  Theorem \ref{thm.boundC}.

\subsection{The bound for degrees of $H_{K,n}$}
Note that
$$
\deg_\pm (fg) = \deg_\pm(f) + \deg_\pm(g),
\quad  \text{and} \quad \maxdeg(f+g) \le \max(\maxdeg(f), \maxdeg(g)).
$$
From $\deg_\pm \{k\} = \pm k/2$, we get
$$
\deg_\pm  (\{k\}!) = \pm k(k+1)/4,
\qquad \deg_\pm  (\qbinom{n}{k}) = \pm k(n-k)/2.
$$
From these and Equation \eqref{33} we get
$$
\maxdeg(H_{K,n}(q)) \le \max_{1\le k \le n+1 }
\left(  -\frac{(2n+2)(2n+3)}{4} + k +
\frac{k}{2}+ \frac{(n+1+k)(n+1-k)}{2} + \maxdeg (J_{K,k}) \right).
$$

Using Proposition \ref{21} for the upper bound of $\maxdeg
(J_{K,k})$, after a simplification, we get
$$
\maxdeg(H_{K,n}(q)) \le \max_{1\le k \le n+1 }
\left(  -\frac{n(n+3)}{2} + \frac{c (k-1)^2}{4} +
\frac{(k-1)s_+}{2} + \frac{|\w|(k^2-1)}{4}\right).
$$
The right hand side reaches maximum when $k=n+1$. Using $\w \le
c+2$,  we have
$$
\maxdeg(H_{K,n}(q)) \le n^2c/2 + n(s_+ +c-1)/2.
$$
A similar calculation shows that
$$
\mindeg(H_{K,n}(q)) \ge -( n^2c/2 + n(s_- +c-1)/2).
$$
If we choose $A_0$ bigger than  $c$ and $|s_\pm +c -1|$, then we
have $|\deg_\pm (H_{K,n}) |\le A_0 n^2$. This proves the first
statement of Theorem \ref{thm.boundC}.

\subsection{The bound for the $l^1$-norm of $H_{K,n}$}

 Multiply both sides of \eqref{33}  by $\{2n+2\}!$,
then use  \eqref{31} and Theorem \ref{95}, we see that there is a constant
$A_7$ such that

\begin{equation}
\lbl{eq.L1C} ||\{2n+2\}!H_{K,n}(q)||_1 \leq  (A_7)^n.
\end{equation}
The polynomials
$$
\tH_{K,n}(q) := q^{A_0n^2} H_{K,n}(q) \qquad
\text{and} \qquad  g(q):=\tH_{K,n}(q) \prod_{j=1}^{2n+2}(1-q^j)
$$
have only non-negative degrees in $q$, with $\maxdeg(\tH_{K,n}) \le
2A_0n^2$:

\begin{equation}
\tH_{K,n}(q) = \sum_{k=0}^{2A_0n^2} a_k q^k
\end{equation}
 Since $g(q)$ is the product of the polynomial on the left hand
side of \eqref{eq.L1C} and a power of $q$, we have

\begin{equation}
 ||g(q)||_1 \leq (A_7)^n.
\label{41}
\end{equation}
There are estimates of $l^1$-norm using Mahler measure \cite{Mahler}.
However, the estimate \eqref{41} is weak: the  inequalities of Mahler
imply an exponential upper bound on the Mahler measure of $H_{K,n}(q)$,
and a doubly exponential upper bound on the $l^1$-norm of $H_{K,n}(q)$. The
following estimate, which does not follow from Mahler
measure considerations, was communicated to us by D. Boyd. Since

$$
\tH_{K,n}(q) = g(q) \frac{1}{\prod_{k=j}^{2n+2} (1-q^j)},
$$
we have that
$$
a_k=\sum_{i=0}^k b_i c_{k-i}, \quad \text{where} \quad g(q) =
\sum_{k} b_k q^k \quad \text{and} \quad \frac{1}{\prod_{k=1}^{2n+2}
(1-q^k)}=\sum_{k=0}^\infty c_k q^k.
$$
Note that $c_k$ is the number of partitions of $k$ of length $\le
2n+2$. Hence $0\le c_{k-1} \le c_k$, and $c_k\le p_k$, where $p_k$
is the number of partitions of $k$. Using the growth rate of $p_k$
(see \cite{An}), we see that there is a constant $A_8$ such that
\begin{equation}
p_k < (A_8)^{\sqrt k}.\label{42}
\end{equation}
The crucial part of the above inequality is the exponent
$\sqrt{k}$. Now we can easily
obtain the desired upper bounds for $||\tH_{K,n}||_1$.
Since $ a_k=\sum_{i=0}^k b_i c_{k-i}$  we have

$$
\begin{aligned}
|a_k| & \le \sum_{i=0}^k |b_i| c_{k-i} \le \left( \sum_{i=0}^k
|b_i|\right) c_{k} \le ||g(q)||_1 c_k\\
&\le (A_7)^n (A_8)^{n \sqrt {2A_0}}  \qquad \text{by \eqref{41},
\eqref{42}  and $k \le 2A_0 n^2 $}
\end{aligned}
$$
It follows that, for $n \ge 1$,

$$
||\tH_{K,n}||_1 \le \sum_{k=0}^{2A_0n^2} |a_k|
\le 2A_0n^2 (A_7)^n (A_8)^{n \sqrt
{2A_0}} \le (A_1)^n,
$$
for appropriate $A_1$. This completes the
proof of Theorem \ref{thm.boundC}.

\section{Growth rates of $R$-matrices and the Lobachevsky function}
\lbl{sec.hypgeom}

\subsection{The Lobachevsky function}
\lbl{sub.estimates}

In Section \ref{sec.boundsJones} we got a simple but crude estimate
for the $l^1$-norm of the $R$-matrices, which  are  a ratio of five
quantum factorials. In this largely independent section we will give
 refined (and optimal) estimates for the growth rate of the
$R$-matrices. These estimates reveal the close relationship between
hyperbolic geometry and the asymptotics of the quantum factorials.

Recall that  the {\em Lobachevsky function} is given by
$$
\La(z)=-\int_0^z \log|2 \sin x| dx = \frac{1}{2} \sum_{n=1}^\infty
\frac{\sin(2 n z)}{n^2}
$$
The Lobachevsky function is an odd, periodic function with period
$\pi$. Its graph for $z \in [0, \pi]$ is:
$$
\psdraw{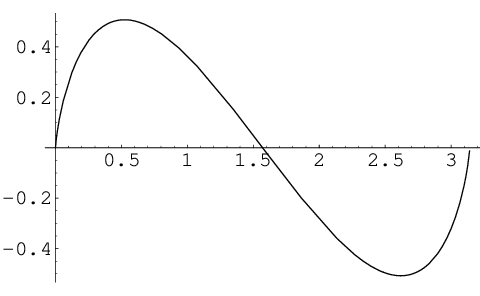}{2.5in}
$$

\begin{definition}
\lbl{eq.ev}
If $f(q) \in \BZ[q^{\pm 1/4}]$, let us denote by $\ev_n(f)$
the {\em evaluation} of $f$ at $q^{1/4}=e^{\pi i/(2n)}$.
\end{definition}
For $0\le k\le n$ we have
$$
\ev_n|\{k\}| =|e^{k\pi i/n}   -  e^{-k\pi i/n}| = 2\sin (k \pi/n),
$$
hence,
$$
\log (\ev_n|(\{j\}!|)|) = \sum_{k=1}^{j} \log|2\sin (k \pi/n)|,
$$
which is very closely related to a Riemann sum of the integral in
the definition of the Lobachevsky function. It is not surprising to
have the following.

\begin{proposition}
\lbl{lem.1}
For every $\a \in (0,1)$ we have:
$$
\log\left|\ev_n (\{ \lfloor \a n \rfloor \}!)\right|=
-\frac{n}{\pi} \La(\pi \a) +O(\log n).
$$
Here $O(\log n)$ is a term which is bounded by $C \log n$ for some
constant $C$ {\em independent} of $\a$.
\end{proposition}

\begin{remark}
\lbl{rem.subleading}
The proof reveals an asymptotic expansion of the form:
$$
\ev_n (\{ \lfloor \a n \rfloor \}!) \sim n^{\theta}
\exp\left(-\frac{n}{\pi} \La(\pi \a) \right)\left(C_0 + \frac{C_1}{n} +
\frac{C_2}{n^2} + \dots\right)
$$
for explicitly computable constants $C_i$ and $\theta$.
\end{remark}

\begin{proof}
Recall the {\em Euler-MacLaurin summation formula},
with error term (see for example, \cite[Chpt. 8]{O}):
$$
\sum_{k=a}^b f(k)=\int_a^b f(x) dx + \frac{1}{2} f(a) +
\frac{1}{2} f(b) +
\sum_{k=1}^{m-1} \frac{B_{2k}}{(2k)!}( f^{(2k-1)}(b)-f^{(2k-1)}(a) )
+ R_m(a,b,f)
$$
where $B_k$ is the $k$th {\em Bernoulli number} and the error term has an
estimate
$$
|R_m(a,b,f)| \leq (2-2^{1-2m}) \frac{|B_{2m}|}{(2m)!}
\int_a^b |f^{(2m)}(x)| dx.
$$

Applying the above formula  for $m=1$ to $f(x)=\log(2 \sin
x\pi/n)$, we have:
\begin{eqnarray*}
\log(\prod_{k=1}^{\lfloor \a n \rfloor } 2\sin (k \pi/n)) & = &
\frac{1}{2} (f(1)+f(\lfloor \a n \rfloor ))+
\int_1^{\lfloor \a n \rfloor } \log(2\sin (t \pi/n)) dt
+R_1(1,\lfloor \a n \rfloor,f) \\
& = &  \frac{1}{2} (f(1)+f( \a n )+
\int_1^{\a n } \log(2\sin (t \pi/n)) dt
+R_1(1,\lfloor \a n \rfloor,f) + \e(\a,n) \\
& = & \frac{1}{2} (f(1)+f(\a n))+ \frac{n}{\pi}  \int_{\pi/n}^{\pi
\a} \log|2 \sin (  u)| u +R_1(1,\lfloor \a n \rfloor, f) +\e(\a,n)
\\
& = & \frac{1}{2} (f(1)+f(\a n))+ \frac{n}{\pi} \left(-\La(\pi \a) +
\La(\frac{\pi}{n}) \right) + R_1(1,\lfloor \a n \rfloor, f) +
\e(\a,n).
\end{eqnarray*}
Here $\e(\a,n)$ comes from adjusting the boundary of integration and
satisfies $|\e(\a,n) | = O(1)$. Note that
$$
\frac{1}{2} |f(1)+f(\a n)| = O( \log n).
$$
Moreover,
$f''(x)=\frac{\pi^2}{n^2} (\csc(\pi x/n))^2 >0$. Hence

$$
\int_1^{\lfloor \al n\rfloor} |f''(x)| dx
=\int_1^{\lfloor \al n\rfloor} f''(x) dx \le \int_1^{\al n} f''(x) dx
= \frac{\pi}{n} \left( \cot(\al \pi) - \cot\left(\frac{\pi}{n}\right)\right).
$$

It follows easily that
$$
|R_1(1,\lfloor \a n \rfloor, f)| = O(1).
$$
Furthermore, using L'Hospital's rule,
one can see  that
$$
\frac{n}{\pi} \left|\La(\frac{\pi}{n})\right| = O( \log n).
$$
The result follows.
\end{proof}

\begin{corollary} For every $\a \in (0,1)$  and any fixed number $d$ we have:
$$
\log\left|\ev_n (\{ \lfloor \a n  +d \rfloor \}!)\right|=
-\frac{n}{\pi} \La(\pi \a) +O(\log n).
$$
\label{87}
\end{corollary}

\begin{proof}
There is $\ve >0$ such that for big enough $n$, we
have $\ve \le x/n \le 1-\ve$ for every integer $x$ between $\lfloor
\a n \rfloor$ and $\lfloor \a n  +d \rfloor$. For such $x$, we have
$ 0< 2 \sin \ve \pi   < 2\sin (x\pi/n)  < 2$, and hence there is a
constant $M$ such that $|\log 2\sin (x\pi/n) | < M$. There are at most
$|d|+1$ such values of $x$. Hence the difference  between
$\log\left|\ev_n (\{ \lfloor \a n  +d \rfloor \}!)\right|$ and $
\log\left|\ev_n (\{ \lfloor \a n \rfloor \}!)\right|$ by absolute
value is less than $(|d|+1) M$, a constant. The result follows.
\end{proof}

\subsection{Asymptotics of the $R$-matrix using ideal octahedra}
\lbl{sub.assRv8}

Since the entries of the $R$-matrix are given by ratios of five quantum
factorials (see Equations \eqref{eq.R+-}), Proposition \ref{lem.1}
gives a formula for the asymptotic behavior of the entries of the
$R$-matrix when evaluated at $e^{2 \pi i/n}$.
This is the content of the next proposition.

\begin{proposition}
\lbl{prop.R+limit}
{\rm (a)}
Suppose that $\a,\b,\k$ are real numbers that satisfy the inequalities
\begin{equation}
\lbl{eq.abk}
\a,\b,\k \in [0,1] \qquad
0 \leq \b+\k \leq 1, \qquad 0 \leq \a-\k \leq 1.
\end{equation}
Then the following limit exists
\begin{equation}
\lbl{eq.r+}
r_+(\a,\b,\k) :=
\lim_{n\to \infty}
\frac{1}{n} \, \log |\ev_n(R_+(n; \lfloor n\a \rfloor,
\lfloor n\b \rfloor,\lfloor n\k \rfloor ))|,
\end{equation}
and is equal to
\begin{equation}
\lbl{eq.r+equal}
r_+(\a,\b,\k)= [-\La(\pi(\b+\k)) +
\La(\pi\b)+\La(\pi\k)-\La(\pi\a)+\La(\pi(\a-\k))]/\pi.
\end{equation}
{\rm (b)} $r_+(\a,\b,\k)$ equals to $1/(2\pi)$ times the volume
of an ideal octahedron with vertices
\begin{equation}
\lbl{eq.vert8}
(0,1,\infty,z_{\kappa},(z_{\b}z_{\kappa}-1)/(z_{\b}-1),z_{\a}) \in
(\BC\setminus\{0,1\})^6
\end{equation}
where
$(z_{\a},z_{\b},z_{\kappa})=(e^{2 \pi i \a},e^{2 \pi i \b}, e^{2 \pi i \kappa})$.
\newline
{\rm (c)} Suppose that $\a,\b,\k$ are real numbers that satisfy:
\begin{equation*}
 \a,\b,\k \in [0,1] \qquad 0 \leq \a+\k \leq 1, \qquad 0
\leq \b-\k \leq 1.
\end{equation*}
Then the following limit exists
$$
r_-(\a,\b,\k) := \lim_{n\to \infty}
\frac{1}{n}\, \log |\ev_n(R_-(n; \lfloor n\a \rfloor,
\lfloor n\b \rfloor,\lfloor n\k \rfloor ))|.
$$
and is equal to
\begin{equation}
\lbl{62}
r_-(\a,\b,\k) = r_+(\b,\a,\k).
\end{equation}
\end{proposition}

\begin{proof}
(a) Observe that
$$
|\ev_n(\{j\})| = |\ev_n(\{n-j\})|= 2 \sin (j \pi/n), \quad
\text{and } |\ev_n(\{n-1\}!)| =  \prod_{j=1}^{n-1} 2\sin (j\pi/n)
=n.
$$

From these, we have that

\begin{equation}
\lbl{61}
|\ev_n(\{j\}!)| = \frac{n}{|\ev_n(\{n-1-j\}!)|}. 
\end{equation}

Using \eqref{eq.R+-} and then \eqref{61}, we have

\begin{equation}\begin{aligned}
 |\ev_n(R_+(n;a,b,k ))|
&= \frac{|\ev_n(\{b+k\}!)|\, |\ev_n(\{n-1+k-a\}!)|  }{ |\ev_n(\{b\}!)|\,
 |\ev_n(\{k\}!)|\,
 |\ev_n(\{n-1-a\}!)|}\\
&=\frac{|\ev_n(\{b+k\}!)| \,|\ev_n(\{a\}!)|  }{ |\ev_n(\{b\}!)|\,
|\ev_n(\{k\}!)|\,
 |\ev_n(\{a-k\}!)|}
 \end{aligned}
 \label{90}
 \end{equation}

Proposition \ref{lem.1} concludes the proof of (a).
(b) was pointed out to us by D. Thurston. Although this fact is not used
in the proof of Proposition \ref{prop.R} nor in the proof of Theorem
\ref{thm.4}, it is an interesting geometric fact. To prove it, recall that
the boundary of 3-dimensional hyperbolic space is $\BC \cup\{\infty\}$.
Let $T_z$ denote the {\em regular ideal tetrahedron} of shape
$z \in \BC-\{0,1\}$. $T_z$ is isometric to the ideal tetrahedron with
ordered vertices
at $0,1,\infty$ and $z$ in the boundary of 3-dimensional hyperbolic space.
For $z \in \BC\setminus\{0,1\}$, the ideal octahedron $T_z$ is isometric
to $T_{1/(1-z)}$ and $T_{z/(z-1)}$ by an orientation-preserving isometry,
and isometric to $T_{1/z}$, $T_{(1-z)/z}$ and $T_{(z-1)/z}$ by an orientation-reversing
isometry. Thus, when $z  \in \BC\setminus\{0,1\}$, we have:
\begin{equation}
\lbl{eq.Tzorbit}
\vol(T_z)=\vol(T_{1/(1-z)})=\vol(T_{z/(z-1)})=-\vol(T_{1/z})=-\vol(T_{(1-z)/z})=
-\vol(T_{(z-1)/z})
\end{equation}
The shape of the ideal tetrahedron
with distinct ordered vertices $(z_0,z_1,z_2,z_3)$ in $\BC\cup\{\infty\}$
is given by the {\em cross-ratio}
\begin{equation}
\lbl{eq.CR}
[z_0:z_1:z_2:z_3]=\frac{(z_0-z_3)(z_1-z_2)}{(z_0-z_2)(z_1-z_3)}
\end{equation}
following the convention of \cite[Eqn.1.4]{DZ}.
If $(\a_1,\a_2,\a_3)$ denote the three dihedral angles of $T_z$ at opposite pairs
of edges, then the volume $\vol(T_z)$ is given by \cite[Thm.10.4.10]{Ra}
$$
\vol(T_z)=\La(\a_1)+\La(\a_2)+\La(\a_3)
$$
If the shape parameter $z=e^{i \th}$ is a complex number of magnitude $1$,
then the dihedral angles of $T_z$ coincide with the angles
of an isosceles triangle with angles $(\th,(\pi-\th)/2,(\pi-\th)/2)$.
In that case, we have
$$
\vol(T_{e^{i \th}})=\La\left(\th\right)+
\La\left(\frac{\pi-\th}{2}\right)+\La\left(\frac{\pi-\th}{2}\right).
$$
The following symmetries of the Lobachevsky function
\cite[Thm.10.4.3,10.4.4]{Ra}
$$
\La(-\th)=-\La(\th), \qquad \frac{1}{2} \La(2\th)=\La(\th)
+\La\left(\th+\frac{\pi}{2}\right)
$$
imply that
\begin{equation}
\lbl{eq.Tzth}
\vol(T_{e^{i \th}})=2 \La\left(\frac{\th}{2}\right).
\end{equation}
Now, we return to the proof of part (b). Fix $\a,\b,\kappa$ as in
\eqref{eq.abk} and consider the complex numbers of magnitude $1$:
$$
(z_{\a},z_{\b},z_{\kappa})=(e^{2 \pi i \a},e^{2 \pi i \b}, e^{2 \pi i \kappa})
$$
Consider five ideal tetrahedra with shapes
\begin{equation}
\lbl{eq.5shapes}
(z_{\b}z_{\kappa})^{-1}, \quad z_{\b}, \quad z_{\kappa}, \quad
z_{\a}^{-1}, \quad z_{\a}z_{\kappa}^{-1}
\end{equation}
Equations \eqref{eq.r+equal} and \eqref{eq.Tzth} implies that the sum of
their volumes is given by $2 \pi \, r_+(\a,\b,\kappa)$. Now consider the
ideal octahedron with vertices
$$
(A,B,C,D,E,F)=(0,1,\infty,z_{\kappa},(z_{\b}z_{\kappa}-1)/(z_{\b}-1),z_{\a})
$$
drawn as follows:
$$
\psdraw{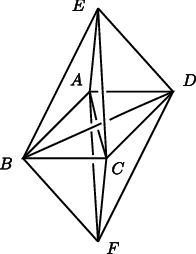}{1in}
$$
It can be triangulated into five ideal tetrahedra
$ABDE$, $BDCE$, $ABCD$, $ABCF$ and $ACDF$ with ordered vertices and
with shape parameters
$$
\left\{1 - \frac{1}{z_{\b} z_{\kappa}}, \frac{z_{\b}}{z_{\b}-1}, \frac{z_{\kappa}}{z_{\kappa}-1}, \frac{z_{\a}}{z_{\a}-1},
\frac{z_{\a}}{z_{\kappa}} \right\}
$$
computed according to Equation \eqref{eq.CR}. Adding up the volumes
of these tetrahedra, with proper orientations concludes the proof of (b).
(c) is analogous to (a). We thank the referee for correcting the vertices
of the octahedron in an earlier version of this paper.
\end{proof}

\subsection{The maximum of the growth rate of the $R$-matrix}

In this section we determine the maximum of $r_\pm(\a,\b,\k)$.

\begin{proposition}
\lbl{prop.R}
{\rm (a)}
With $\a,\b,\k$ satisfying \eqref{eq.abk}, $r_+(\a,\b,\k)$ achieves
maximum when $\a=3/4, \b=1/4$, and $\k=1/2$. Moreover
$$
r_+(3/4,1/4,1/2)=\frac{v_8}{2 \pi},
$$
where
$$
v_8=8 \La(\pi/4) \approx 3.6638623767088760602 \dots
$$
is the volume of the regular hyperbolic ideal octahedron.
\newline
{\rm (b)}
Similarly,  $r_-(\a,\b,\k)$ reaches maximum when $\a=1/4, \b=3/4$,
and $k=1/2$; and its maximum value is the same as that of
$r_+(\a,\b,\k)$.
\end{proposition}

Thus, asymptotically, the winning configuration is given by:
$$
\psdraw{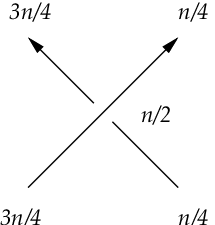}{1in}
$$

\begin{proof}
It is enough to consider the case of $r_+$. The result
for $r_-$ follows from \eqref{62}. Let $\d = \a -\k$, we have
$$
r_+(\a,\b,\k) = -\La(\pi(\b+\k)) +
\La(\pi\b)+\La(\pi\k)-\La(\pi(\d+ \k))+\La(\pi(\d)),
$$
with domain $ 0\le \b,\d,\k$, and $ \b+\k \le 1, \d+\k \le 1$. Note
the symmetry between $\b$ and $\d$.

Using $\La'(x) = - \log(2\sin x)$ for $ 0 <x < \pi$, one can easily
show that the function $-\La(\pi(\b+\k)) + \La(\pi\b)$, for a fixed
$\k \in [0,1]$, achieves maximum at $ \b = (1-\k)/2$. It follows
that the maximum of $r_+(\a,\b,\k)$ is the same as the maximum of

$$
g(\k):= 2(-\La(\pi(\b+\k)) + \La(\pi\b)) + \La(\pi\k),
$$
with $ \b = (1-\k)/2$. The domain for $g$ is $\k \in [0,1]$. Using
the derivative of $g$ it is easy to show that $g$ achieves maximum
when $\k=1/2$. In this case $\a=3/4, \b=1/4$.
\end{proof}

\begin{remark}
Another proof is to use part (b) of Proposition \ref{prop.R+limit} and
the fact that the volume of an ideal octahedron is maximized at a regular
ideal octahedron; see \cite{Ra}.
\end{remark}

\begin{lemma}
\lbl{lem.2circles}
If $z,w$ are complex numbers that satisfy $|z|=|w|=1$ and $|1-z|=|1-w|$,
then $z=w^{\pm 1}$.
\end{lemma}

\begin{proof}
Let us define
$$
C_{u_0,r}:=\{ u \in \BC \,\, | \,\, |u-u_0|=r>0 \}.
$$
Then $C_{u_0,r}$ is a circle with center $u_0$ and radius $r$. Fixing $w$, it
follows that $z \in C_{0,1} \cap C_{1,|1-w|}$. The intersection of two circles
is two points, and since $w$ and $w^{-1}=\bar{w}$ both lie in the intersection,
the result follows.
\end{proof}

\subsection{The maximum of the $R$-matrices at roots of unity}

Proposition \ref{prop.R} gives the maximum of the growth rate of
$\ev_n(R_+(n;a,b,k))$, as $n\to \infty$. The following proposition
gives the maximum of $\ev_n(R_+(n;a,b,k))$, for a fixed $n$.

\begin{proposition}
\lbl{prop.Rpm} {\rm (Proof in Section \ref{sub.pf71})}
The value of
$|\ev_n(R_+(n;a,b,k))|$ achieves maximum at $a= \lfloor 3n/4
\rfloor$, $b= \lfloor (n-1)/4 \rfloor$, and $k= a-b$. The value of
$|\ev_n(R_-(n;a,b,k))|$ achieves maximum at $a= \lfloor (n-1)/4
\rfloor$, $b= \lfloor 3n/4 \rfloor$, and $k= b-a$. The maximum
value of $|\ev_n(R_+(n;a,b,k))|$ is the same as that of
$|\ev_n(R_-(n;a,b,k))|$. 
\end{proposition}

Note that for these optimal values in the $R_+$ case, $ |a-3n/4| \le
1$, $|b-n/4| \le 1$ and $|k-n/2| \le 1$. The proof of this
proposition will be given in Appendix \ref{CCC}.

From Corollary \ref{87} and Propositions \ref{prop.R+limit},
\ref{prop.R} and \ref{prop.Rpm}, we have the following.

\begin{corollary}\label{86}
The growth rate of the maximum of
$|\ev_n(R_+(n;a,b,k))|$ is given by
$$ \lim_{n\to \infty} \frac{\max_{a, b, k} \log |\ev_n(R_\pm(n;a,b,k))| }{n}
= \frac{v_8}{2 \pi}.
$$
\end{corollary}

\subsection{Proof of Theorem \ref{thm.4}}
\lbl{sub.thm4}

Recall that by \eqref{eq.statesum}, the colored Jones function is
the sum of $n^c$ summands.   Each summand  $F(n,\bk)$ is the product
of $R$-matrices (which are weights of crossing points) and weights
of extreme points (which have absolute value 1). There are $c+2$
crossing points, but the weights of the two boundary crossing have
absolute value 1. Hence

$$
| J_{K,n}(e^{2\pi i/n})| \le n^c \left (\max_{a, b, k}
|\ev_n(R_\pm(n;a,b,k))|\right)^{c}.
$$

From the growth rate of $\max_{a, b, k} |\ev_n(R_\pm(n;a,b,k))|$
given by Corollary \ref{86} we get the theorem.

\section{The generalized volume conjecture near $\a=2\pi i$}
\lbl{sec.cyclotomic}

In this section we will prove Theorems \ref{thm.near1} and
\ref{thm.near2} which are concerned with the Generalized Volume
Conjecture near $2 \pi i$. Our proofs use crucially  the well-known
{\em symmetry principle}, see \cite{KM,Le1}: Suppose  $m, m'$ and
$n$ are positive integers with $m \equiv \pm m' \bmod n$, then
\begin{equation}
\lbl{eq.symmetryJ} J_{K,m}(e^{2 \pi i/n})=J_{K,m'}(e^{2 \pi i/n}).
\end{equation}
Note that this fact  is also a consequence of the existence of the
cyclotomic expansion. However, the case of higher rank Lie algebra
requires results from canonical basis theory, see \cite{Le1}.

\begin{proof}(of Theorem \ref{thm.near1})
The symmetry principle implies that for all $n  >m >0$, we have:
$$
J_{K,n\pm m}(e^{2 \pi i/n})=J_{K,m}(e^{2 \pi i/n})
$$
which implies that
$$
\lim_{n\to\infty}J_{K,n\pm m}(e^{2 \pi i/n})=\lim_{n\to\infty}
J_{K,m}(e^{2 \pi i/n})=J_{K,m}(1)=1,
$$
from which Theorem \ref{thm.near1} follows easily.
\end{proof}

\begin{proof}(of Theorem \ref{thm.near2})
Fix a knot $K$ and consider the neighborhood $U_K$ of $0$ as in
Theorem \ref{thm.11}. Define $V_K=1+U_K$.

Let us suppose that $\a/(2 \pi i) \in V_K$ is a rational number not
equal to $1$. Assume that  $\a=2 \pi i p/m$ with $p,m$ unequal
coprime positive integers.  Let $N=np$. Then, the symmetry principle
implies that
\begin{eqnarray*}
f_{K,N}(\a) &=& J_{K,N}(e^{\a/N}) \\
&=& J_{K,np}(e^{2 \pi i/(nm)}) \\
&=& J_{K, n|p-m|}(e^{2 \pi i/(nm)}).
\end{eqnarray*}
Since $n|p-m|/(nm)=|p/m-1|\in U_K$, Theorem \ref{thm.11} implies
that
$$
\lim_{n\to\infty} J_{K, n|p-m|}(e^{2 \pi i/(nm)}) = \frac{1}{\D(e^{2
\pi i(|p/m-1|})}.
$$
In other words,
$$
\lim_{n\to\infty} f_{K,np}(\a)= \frac{1}{\D(e^{2 \pi i(|p/m-1|})}$$
is bounded. The result follows.
\end{proof}

\section{The $q$-holonomic point of view}
\lbl{sec.qholonomic}

\subsection{Bounds on $l^1$-norm  of $q$-holonomic functions}
\lbl{sub.qholo}

The main result of \cite{GL} is that for every knot $K$, the
functions $J_K$ and $H_K$ are $q$-holonomic. Recall that a sequence
$f:\BN\longto\BQ(q)$ is $q$-{\em holonomic} if satisfies a
$q$-{linear difference equation}. In other words, there exists a
natural number $d$ and polynomial $a_j(u,v) \in \BQ[u,v]$ for
$j=0,\dots,d$ with $a_d \neq 0$ such that for all $n \in \BN$ we
have:
\begin{equation}
\sum_{j=0}^d a_j(q^n,q) f_{n+j}(q)=0.
\label{e0001}
\end{equation}

In this section we observe that $q$-holonomic functions satisfy {\em
a priori} upper bounds on their degrees and (under an integrality
assumption) on their $l^1$-norm. As a simple corollary, we obtain
another proof of the  quadratic bounds in   Proposition \ref{21}, though not as explicit.

\begin{definition}
\lbl{eq.integral} We say that a sequence $f:\BN\longto\BZ[q^{\pm 1}]$
is $q$-{\em integral holonomic} if it satisfies an {\em
$q$-difference equation} as above with $a_d=1$.
\end{definition}

\begin{question}
\lbl{que.1}
Is it true that $J_K$ and $C_K$ are $q$-integral holonomic for every knot $K$?
\end{question}

For a partial answer, see \cite{GS}.

\begin{theorem}
\lbl{thm.L1bound}
$\mathrm{(a)}$
If $f:\BN\longto \BZ[q^{\pm 1}]$ is $q$-holonomic, then for all $n$ we have:
$$
\maxdeg(f_n)=O(n^2)   \qquad
\text{and} \qquad \mindeg(f_n) =O(n^2).
$$
\newline
$\mathrm{(b)}$ If $f$ is $q$-integral holonomic, then for all $n$ we have:
$$
||f_n||_1 \leq C^n
$$
for some constant $C$. In particular,

$$
\limsup_{n \to \infty} \frac{\log|f_n(e^{\a/n})|}{n}
\leq C_{\a}
$$
for all $\a \in \BC$.
\end{theorem}

In other words, integral $q$-holonomic functions grow at most
exponentially.

\begin{proof}
Suppose $f$ satisfies \eqref{e0001}.
It is easy to see that for every $a(u,v)\in \BQ[u,v]$, there exists a 
constant $C'$ such that $\deg_+(a(q^n,q)) < C' n $  for every $n \ge 1$. 
We choose such a common $C'$ for  all
$a_j(q^n,q)$, $j=0,1,\dots,d$, and, in addition,  $C' > \maxdeg f(n)$ for 
$n=0, 1,\dots,d$.

We will prove by induction on $n\ge 1 $ that $\maxdeg f(n) \leq C' n^2$.
By assumption, it is true for $n=1, \dots, d$. For $n \ge 1$,
Then, by induction we have:
\begin{eqnarray*}
\maxdeg  f_{n+d}(q) & = &
\maxdeg \Big(  a_d(q^n,q) \, f_{n+d}(q)\Big) - \maxdeg  a_d(q^n,q) \\
&=& \maxdeg  \left(-\sum_{j=0}^{d-1} a_j(q^n,q) f_{n+j}(q)\right)   - \maxdeg  
a_d(q^n,q) \\
&< & C'n + C'(n+d-1)^2 + C'n \le C'(n+d)^2.
\end{eqnarray*}
The second claim in (a) follows similarly.

For (b), let $c_j=||a_j(Q,q)||_1$ for $j=0, \dots, d-1$, and choose $C$
so that
\begin{itemize}
\item
$C^d \geq c_{d-1} C^{d-1} + \dots + c_0 C^0$, and
\item
$||f_n(q)||_1 \leq C^n$ for $n=0, \dots, d-1$.
\end{itemize}
Then, it is easy to see by induction that (b) holds for all $n$.
\end{proof}

\begin{remark}
\lbl{rem.holosharp}
It is easy to see that the bounds of Theorem \ref{thm.L1bound} are sharp.
For example,
consider the sequence $f_n(q)=(1+q)(1+q^2)\dots(1+q^n)$.
\end{remark}

 Theorem \ref{thm.L1bound} gives an alternative proof of the quadratic bounds for the degrees of the color Jones polynomial, though not as
 explicit as in
Proposition \ref{21}.  If Question
\ref{que.1} has a positive answer then Theorem \ref{thm.L1bound} also gives an alternative proof of Theorem \ref{thm.3}.


\subsection{Bounds for higher rank groups}
\lbl{sub.high}

In \cite{GL}, we considered the colored Jones function
$$
J_{\fg,K}: \La_w \longto \BZ[q^{\pm 1}]
$$
of a knot $K$, where $\fg$ is a {\em simple Lie algebra} with {\em
weight lattice}
$\La_w$. In the above reference, the authors proved that $J_{\fg,K}$
is a $q$-holonomic function, at least when $\fg$ is not $G_2$.
For $\fg=\mathfrak{sl}_2$, $J_{\mathfrak{sl}_2,K}$ is the colored Jones
function $J_K$ discussed earlier.

In \cite{GL} , the authors gave  state-sum formulas for $J_{\fg,K}$
similar to \eqref{eq.statesum} where the summand takes values in
$\BZ[q^{\pm 1/D}]$, where $D$ is the size of the center of $\fg$.

The methods of the present paper give an upper bound for
the growth-rate of the $\fg$-colored Jones function. More precisely,
we have:

\begin{theorem}
\lbl{thm.gJones}
For every simple Lie algebra $\fg$ (other than $G_2$)  and every
$\a \in \BC$, and every $\l \in \La_w$, there exists a
constant $C_{\fg,\a,\l}$ such that for every knot with $c+2$ crossings,
we have:
$$
\limsup_{n \to \infty} \frac{\log|J_{\fg,K,n\l}(e^{\a/n})|}{n} \leq
C_{\fg,\a,\l} c .
$$
\end{theorem}

The details of the above theorem will be explained in a subsequent
publication.

\section{Some physics}
\lbl{sec.physics}

\subsection{A small dose of physics}
\lbl{sub.physics}

One does not need to know the relation of the colored Jones
function and quantum field theory in order to understand the statement
and proof of Theorem \ref{thm.1}. Nevertheless, we want to add some
philosophical comments, for the benefit of the willing reader.
According to Witten (see \cite{Wi}), the Jones polynomial
$J_{K,n}$ can be expressed by a partition function of a topological quantum
field theory in $3$ dimensions---a gauge theory with Chern-Simons Lagrangian.
The stationary points of the Lagrangian
correspond to $SU(2)$-flat connections on an ambient manifold, and the
observables are knots, colored by the $n$-dimensional irreducible
representation of $SU(2)$. In case of a knot in $S^3$, there is only one
ambient flat connection, and the corresponding perturbation theory
is a formal power series in $h=\log q$.

Rozansky exploited a cut-and-paste property of the Chern-Simons path integral
and considered perturbation theory of the knot complement, along an abelian
flat connection with monodromy given by \eqref{eq.meridian}. In fact, Rozansky
calls such an expansion the $U(1)$-{\em reducible connection contribution}
(in short, $U(1)$-RCC) to the
Chern-Simons path integral, where $U(1)$ stands for the fact that the flat
$SU(2)$ connections are actually $U(1)$-valued abelian connections.
Formal properties of such a perturbative expansion, enabled Rozansky to
deduce (in physics terms) the loop expansion of the colored Jones function;
see \cite{Ro2}.
Rozansky also proved the existence of the loop expansion
using an explicit state-sum description of the colored Jones function;
see \cite{Ro1}.

Of course, perturbation theory means studying formal power series that rarely
converge. Perturbation theory at the trivial flat connection in a knot
complement converges, as it resums to a Laurent polynomial in $e^h$; namely
the $n$th colored Jones polynomial.
The volume conjecture for small complex angles is precisely the
statement that perturbation theory for abelian flat connections (near the
trivial one) does converge.

At the moment, there is no physics (or otherwise) formulation of perturbation
theory of the Chern-Simons path integral along a discrete and faithful
$\SL_2(\BC)$ representation. Nor is there an adequate explanation of the
relation between $SU(2)$ gauge theory (valid near $\a=0$) and
a complexified $\SL_2(\BC)$ gauge theory, valid near $\a=2\pi i$.
These are important and tantalizing questions, with no answers at present.

\subsection{The WKB method}
\lbl{sec.WKB}

Since we are discussing physics interpretations of Theorem \ref{thm.1}
let us make some more comments. Obviously, when the angle
$\a$ is sufficiently big, the asymptotic expansion of Equation \eqref{eq.thm1}
may break down. For example, when $e^{\a}$ is a complex root of the Alexander
polynomial, then the right hand side of \eqref{eq.thm1} does not make sense,
even to leading order. In fact, when $\a$ is near $2 \pi i$, then the solutions
are expected to grow exponentially, and not polynomially, according to the
Volume Conjecture.

The breakdown and change of rate of asymptotics is
a well-documented phenomenon well-known in physics, associated with WKB
analysis, after Wentzel-Krammer-Brillouin; see for example \cite{O}.
 In fact, one may obtain an independent proof of Theorem \ref{thm.1}
using {\em WKB analysis}, that is, the study of asymptotics of solutions of
difference equations with a small parameter.
The key idea is that the sequence of colored Jones
functions is a solution of a linear $q$-difference equation, as was
established in \cite{GL}. A discussion on WKB analysis of $q$-difference
equations was given by Geronimo and the first author in \cite{GG}.

The WKB analysis can, in particular, determine {\em small exponential
corrections} of the form $e^{-c_{\a} n}$ to the asymptotic expansion of
Theorem \ref{thm.1}, where $c_{\a}$ depends on $\a$, with $\text{Re}(c_{\a})
<0$ for $\a$ sufficiently small. These small exponential corrections
(often associated with instantons)
cannot be captured by classical asymptotic analysis (since they vanish
to all orders in $n$), but they are important and dominant (i.e.,
$\text{Re}(c_{\a})>0)$ when $\a$ is
near $2 \pi i$, according to the volume conjecture. Understanding the change of
sign of $\text{Re}(c_{\a})$ past certain so-called Stokes directions
is an important question that WKB addresses.

We will not elaborate or use the WKB analysis in the present paper.
Let us only mention that the loop expansion of the colored Jones function
can be interpreted as WKB asymptotics on a $q$-difference equation satisfied
by the colored Jones function.

\appendix

\section{The volume conjecture for the Borromean rings}
\lbl{sec.borromean}

It is well-known that the complement of the Borromean rings $B$ can
be geometrically identified by gluing two regular ideal octahedra,
see \cite{Th}. As a result, the volume $\vol(S^3-B)$ of $S^3-B$ is
equal to $2 v_8$.

Suppose $L$ is a $k$-component framed link,  and $n_1,\dots,n_k$ are
positive integers. The colored Jones polynomial $\tilde
J_{L}(n_1,\dots,n_k) \in \BZ[q^{\pm 1/4}]$ is the $sl_2$-quantum
invariant of the link whose components are colored by $sl_2$-modules
of dimensions $n_1,\dots,n_k$, see \cite{RT,Tu}. The normalization
is chosen so that for the unknot, $\tilde J_{L}(n) = [n]$. Define

$$J_{L,n}(q) := \frac{J_L(n,n,\dots,n)}{[n]}.$$

The next theorem confirms the volume conjecture for the Borromean rings.

\begin{theorem}
\lbl{thm.borromean} Let $B$ be the Borromean rings, then
\begin{equation*}
\lim_{n \to \infty} \frac{\log|J_{B,n}(e^{2 \pi i/n})|}{n}
=\frac{1}{2 \pi} \vol(S^3-B).
\end{equation*}
\end{theorem}

\begin{proof}
For  an integer $j$ and a positive integer $k$
let $x_j= 2\sin(j\pi/n)$ and $z_k= \prod_{j=1}^{k} x_j$.

Then, see \eqref{61},

\begin{eqnarray}
x_j &=&x_{n-j} = - x_{n+j},
\label{e007}\\
    z_k &=& n/z_{n-1-k} \quad \text{ for $1\le
k\le n-1$}.
\label{44}\end{eqnarray}

Using  Habiro's formula for $\tilde J_L$ of the Borromean ring
\cite{H2,H1},  one has
$$
J_{B,n}(q)= \sum_{l=0}^{n-1}(-1)^l\, \frac{\{n\}^2 \left(
\prod_{j=1}^{l} \{n+j\}
\{n-j\}\right)^3}{\left(\prod_{j=l+1}^{2l+1}\{j\}\right)^2}.
$$

When $q^{1/2}= e^{i \pi/n}$, one has $\{j\}= 2i \sin \frac{j\pi}{n}$,
which is 0 exactly when $j$ is divisible by $n$. Hence if $2l+1 <n$,
then the denominator of the term in the above sum is never 0, while
the numerator is 0, since it has 2 factors $\{n\}$. On the other
hand, if $2l+1 >n$, then the denominator has 2 factors $\{n\}$,
which would cancel with the 2 same  factors of the numerator. Hence
at $q^{1/2}= e^{i \pi/n}$ one can assume that $2l +1 \ge n$, or $l>
n/2-1$:
\begin{eqnarray*}
J_{B,n}(e^{2 \pi i/n})&=& \sum_{n>l>n/2-1}(-1)^l\,\ev_n \frac{ \left(
\prod_{j=1}^{l} \{n+j\}
\{n-j\}\right)^3}{\left(\prod_{j=l+1}^{n-1}\{j\}\,
\prod_{j=n+1}^{2l+1}\{j\}\right)^2}\\
&=& \sum_{n>l>n/2-1} \frac{ \left(
z_l \right)^6}{( z_{n-l-1})^2\, (z_{2l+1-n})^2}\quad \text{by \eqref{e007}}
\end{eqnarray*}

Using  \eqref{44}, which says  $z_l = n/ z_{n-1-l}$, we have

\begin{equation}
J_{B,n}(e^{2 \pi i/n})= \sum_{n>l>n/2-1} n^2\, (\g_l)^2,
\quad \text{where} \quad \g_l= \frac{ \left(
z_{l}\right)^2}{(z_{n-1-l})^2\, z_{2l+1-n}}.
\label{e201}
\end{equation}


By \eqref{84} below, with $a_l= l$, $b_l= n-1-l$ and $k_l=2l+1-n$, we have:
$$
\gamma_l=|\ev_n(R_+(n;a_l,b_l,k_l))|.
$$
 By Proposition \ref{prop.Rpm}, $|\ev_n(R_+(n;a_l,b_l,k_l))|$ achieves
maximum at
$$
a_{\max{}} = \lfloor 3n/4 \rfloor, \quad b_{\max{}} = \lfloor
(n-1)/4 \rfloor, \quad \text{and } k_{\max{}} = a_{\max{}} -b_{\max{}}.$$

 When $l= \lfloor
3n/4 \rfloor$ we have $a_l= a_{\max{}}$, while
$|b_l - b_{\max{}}| \le 1$ and $| k_l - k_{\max{}}| \le 1$.  It is easy to see that

\begin{equation}
\lim_{n\to \infty} \frac{\log \g_{\lfloor
(n-1)/4 \rfloor} -  \log |\ev_n(R(n;a_{\max{}},b_{\max{}},k_{\max{}}))|}{n}=0.
\label{e203}
\end{equation}

There are less than $n$ summands in the right hand side of \eqref{e201},
and each summand is positive.  Hence
$$
 n^2\,  \left( \g _{\lfloor
3n/4 \rfloor}\right) ^2 < J_{B,n}(e^{2 \pi i/n}) < n^3
 |\ev_n(R(n;a_{\max{}},b_{\max{}},k_{\max{}}))|^2
$$
From \eqref{e203} it follows that
$$
2\pi \lim_{n \to \infty} \frac{\log|J_{B,n}(e^{2 \pi i/n})|}{ n} = 2\pi
\lim_{n \to \infty} \frac{\log |\ev_n(R(n;a_{\max{}},b_{\max{}},k_{\max{}}))|^2}{n},
$$
which is equal to $2v_8$, according to Corollary \ref{86}.
\end{proof}

\section{Proof of Proposition \ref{prop.Rpm}}
\lbl{CCC}

\subsection{Preliminary estimates}
\lbl{sub.prelim71}

Again we denote $x_j= 2\sin (j \pi /n)$. The following is obvious.

\begin{lemma}
\lbl{85}
\rm{(a)}
The $x_j$, as a function of $j$,  is increasing for $j\in
[0,n/2]$ and   decreasing for $j \in [n/2,n]$. In particular, for
$j\le l \le n-j$ we have $x_j \le x_l$.
\newline
\rm{(b)}
For every $1\le j\le n-1$, one has $2 \ge x_j$. For $ n/4\le
j\le 3n/4$, one has $2 \le (x_j)^2$.
\end{lemma}

\begin{lemma}
\lbl{45}
For a fixed $k$, $1\le k\le n-1$, the value of\,  $y_b(k):=
\prod_{j=b+1}^{b+k} x_j$ achieves maximum at 
\begin{equation}
b = \beta(k):=\lfloor (n-k)/2
\rfloor.
\end{equation}
\end{lemma}

\begin{proof} 
We will  prove that if $b < \beta(k)$, then $y_b(k) \le y_{b+1}(k)$, while if $b> \beta(k)$ then $y_b(k) \le y_{b-1}(k)$. This will prove the lemma.

Suppose  $b < \beta(k)$. Then  $b \le  \lfloor
(n-k)/2 \rfloor -1 \le (n-k)/2 -1$. It follows that $(b+1) \le b+k+1
\le n-(b+1).$ From Lemma \ref{85}(a) we get $x_{b+1} \le x_{b+k+1}$.
Hence $y_{b+1}(k)/y_b(k) = x_{b+k+1}/x_{b+1} \ge 1$, or $y_{b+1}(k) \ge y_b(k)$.

Suppose now $b \ge 1+ \beta(k)$. If $b \ge n/2$, then
$x_b \ge x_{b+k}$ since $y_j$ is deceasing on $[n/2,n]$. If $ b
<n/2$, then from $b \ge 1+ \lfloor (n-k)/2 \rfloor$ one can easily
show that $b+k \ge n-b \ge n/2$. Hence we also have $x_{b}=
x_{n-b}\ge x_{b+k}$. Thus $y_{b-1}(k)/y_b(k) = x_b /x_{b+k} \ge 1$.
\end{proof}

 Using \eqref{90}, with
$|\ev_n(\{j\})|=x_j$, we have

\begin{equation}
\lbl{84}
|\ev_n R_+(n;a,b,k)| =   \frac{y_b(k) \, y_{a-k}(k) }{y_0(k)
}, \quad |\ev_n R_-(n;a,b,k)| =   \frac{y_a(k) \, y_{b-k}(k) }{y_0(k)
}.
\end{equation}

By Lemma \ref{45}, both
$y_b(k)$ and $y_{a-k}(k)$ achieve maximum when $b=a-k=\beta(k)$. Hence

\begin{equation}
 \max |\ev_n R_+(n;a,b,k)| = \max_{0\le k \le n} s(k) ,
\quad \text{ where } s(k)= \frac{\left(y_{\beta(k)}(k)\right) ^2}{y_0(k)}.
 \lbl{e103}
 \end{equation}

\begin{lemma}
\label{l108}
One has
\begin{equation} \frac{s(k+1)}{s(k)} = \frac{(x_{\beta(k)})^2}{x_{k+1}}
\label{e108}
 \end{equation}
with the denominator satisfying
\begin{equation}
 x_{k+1} = \begin{cases}
            x_{2\beta(k)-1} \quad &\text{ if  $n-k$ is even}\\
 x_{2\beta(k)} \quad & \text{ if  $n-k$ is odd}
           \end{cases}
           \label{e116}
\end{equation}
\end{lemma}

\begin{proof}
By definition
\begin{equation}
s(k) = \frac{\left(y_{\beta(k)}(k)\right) ^2}{y_0(k)}
= \frac{\left( \prod_{j=1}^k  x_{\beta(k)+j}\right) ^2 }{\prod_{j=1}^k x_j}.
\label{e107}
\end{equation}

Note that, with $\beta(k)= \lfloor (n-k)/2 \rfloor$, we have
\begin{equation} \beta(k+1) = \begin{cases} \beta(k)-1 \quad &\text{if $n-k$ is even} \\
 \beta(k) \quad &\text{if $n-k$ is odd}
\end{cases}
\label{e111}
\end{equation}

We consider two cases: $n-k$ is even and $n-k$ odd.

(a) $n-k$ is even.  Replacing $k$ with $k+1$ in
\eqref{e107}, then using $\beta(k+1)= \beta(k)-1$, we get

$$
s(k+1) = \frac{\left( \prod_{j=1}^{k+1}  x_{\beta(k+1)+j}\right) ^2 }{\prod_{j=1}^{k+1} x_j}
= \frac{\left( \prod_{j=1}^{k+1}  x_{b_{k}-1+j}\right) ^2 }{\prod_{j=1}^{k+1} x_j}
= \frac{(x_{\beta(k)})^2 \, \left( \prod_{j=1}^k  x_{\beta(k)+j}\right) ^2 }{
x_{k+1} \prod_{j=1}^{k} x_j}.
$$
Dividing by $s(k)$, we get \eqref{e108}.
As for the denominator, using $x_j= x_{n-j}$ and $n-k= 2\beta(k)$,
$$ x_{k+1} = {x_{n-k-1}} ={x_{2\beta(k)- 1}}.$$
This proves the lemma when $n-k$ is even.

(b) $n-k$ is odd. 
Replacing $k$ with $k+1$ in \eqref{e107}, then using $\beta(k+1)= \beta(k)$,
we get
\begin{equation}
\lbl{e110}
s(k+1)  = \frac{\left( \prod_{j=1}^{k+1}  x_{b_{k}+j}\right)^2 }{\prod_{j=1}^{k+1} x_j}
= \frac{(x_{\beta(k)+k+1})^2 \, \left( \prod_{j=1}^k  x_{\beta(k)+j}\right) ^2 }{x_{k+1}
\prod_{j=1}^{k} x_j}=  \frac{(x_{\beta(k)+k+1})^2 }{x_{k+1} } s(k).
\end{equation}
Using $x_j = x_{n-j}$ and $n-\beta(k) -k-1 = \beta(k)$,
we have
$$
x_{\beta(k) + k+ 1} = x_{n-\beta(k) -k-1} = x_{\beta(k)},
$$
which, together with Equation \eqref{e110}, proves Equation \eqref{e108}.
As for the denominator, using $n-k-1 = 2\beta(k)$,
$$
x_{k+1} = x_{n-k-1} = x_{2 \beta(k)}.
$$
This completes the proof of the lemma.
\end{proof}

As $k$ increases from $0$ to $n$, $\beta(k) = \lfloor (n-k)/2 \rfloor$ decreases and covers all integers from  $\lfloor n/2 \rfloor$ to 0.  

\begin{lemma}\label{l115} (a) If  $n \ge 7$ then $s(k)$ achieves maximum at an integer $k$ such that $\beta(k) = \lfloor (n-1)/4 \rfloor$.

(b) $s(k)$ achieves maximum at $k$ which is the smallest integer such that $\beta(k) = \lfloor (n-1)/4 \rfloor$.
\end{lemma}

\begin{proof} (a)
We will  show that

(a1) if $\beta(k) >(n-1)/4$ then $s(k+1) \ge  s(k)$, and

(a2) if $ \beta(k) \le  (n-1)/4-1$ then
 $s(k-1) >  s(k)$.
 
This will show that the maximum can be achieved for a $k$ such that  $\beta(k)= \lfloor (n-1)/4 \rfloor$.

Proof of (a1).  Suppose $\beta(k) > \frac{n-1} 4 $. There is no integer in the interval $(\frac{n-1} 4,\frac n 4)$, because otherwise (by multiplying by 4) there would be an integer in $(n-1,n)$. It follows that
$\beta(k) \ge n/4$.

Besides, $\beta(k) = \lfloor (n-k)/2 \rfloor \le n/2$. Thus $\beta(k) \in [n/4, n/2]$.
By Lemma \ref{85}(b), $(x_{\beta(k)})^2 \ge 2 \ge x_{j}$ for any $ 1\le j \le n$.
It follows that the right hand side of \eqref{e108} is bigger than or equal
to 1, or $s(k+1)/s(k) \ge 1$.

Proof of (a2). Suppose $\beta(k) \le  \frac{n-1} 4 -1$. Then  $\beta(k-1) \le  \frac{n-1} 4 $ since by \eqref{e111}, either $\beta(k-1)= \beta(k)$ or $\beta(k-1)=\beta(k)+1$.

Since   $2\beta(k-1) < n/2$, by Lemma \ref{85}(a), $x_{2\beta(k-1)} \ge x_{2\beta(k-1)-1}$.
It follows that $x_k$, either equal to   $x_{2\beta(k-1)}$ or $x_{2\beta(k-1)-1}$ by \eqref{e116}, satisfies
\begin{equation*}
 x_k \ge x_{2\beta(k-1)-1}.
\end{equation*}
By Lemma \ref{l108} and the above inequality, 

$$
\frac{s(k)}{s(k-1)} = \frac{(x_{\beta(k-1)})^2}{x_{k}}
\le \frac{(x_{\beta(k-1)})^2}{x_{2\beta(k-1)-1}} < 1
$$
where the last inequality follows from Lemma \ref{leml101} below.
This completes the proof of the part (a). 

(b) When $n < 7$ the statement is checked by explicit calculation. We will assume $n\ge 7$.

There are two values of $k$ such that  $\beta(k)= \lfloor (n-1)/4 \rfloor$. Let $\kappa$ be the smaller one, then the other one is $\kappa+1$.
Then $n-\kappa$ is odd since otherwise $\beta(\kappa-1)=\beta(k)$. 

Since $n-\kappa$ is odd, by Lemma  \ref{l108}, we have 

$$ \frac{s(\kappa+1)}{s(\kappa)} = \frac{(x_{\beta(\kappa)})^2}{x_{2\beta(\kappa)-1}},$$
which is less than 1 by Lemma \ref{leml101}. This means $s(k)$ achieves maximum at $k=\kappa$.
\end{proof}

\subsection{Proof of  Proposition \ref{prop.Rpm}}
\lbl{sub.pf71}
First consider $R_+(n;a,b,k)$. By \eqref{e103}, Lemmas \ref{l115}(b) and \ref{45}, $|R_+(n;a,b,k)|$ achieves maximum when $k$ satisfies the condition  in Lemma \ref{l115}(b),
$b = \lfloor (n-1)/4 \rfloor$, and $a= b+k$. The value of $k$ satisfying the condition in Lemma \ref{l115}(b) can be calculated easily:

$$
k= \begin{cases}
n/2 + 1  &   n=0 \bmod 4 \\
(n-1)/2  &  n=1 \bmod 4 \\
n/2      &  n=2 \bmod 4 \\
(n+1)/2  &  n=3 \bmod 4
\end{cases}
$$
From there one can calculate $a= k+b$. It is easy to check that the values of $a,b,k$ are exactly the ones given in Proposition \ref{prop.Rpm}. 

Now turn to 
 $R_-(n;a,b,k)$. By \eqref{84}, 
 $$ |R_-(n;a,b,k)| = |R_+(n;b,a,k)|.$$
Hence $|R_-(n;a,b,k)|$  and $|R_+(n;b,a,k)|$ have the same maximum, and $|R_-(n;a,b,k)|$ achieves maximum when $|R_+(n;b,a,k)|$ achieves maximum.

This completes the proof of
Proposition \ref{prop.Rpm}, modulo the following lemma.

\begin{lemma}
\lbl{leml101}
For $1 \le j \le \frac {n-1}4$, with $n \ge 7$, one has $x_{2j-1} > x_j^2$.
\end{lemma}

\begin{proof} With $x_j = 2 \sin (j \pi/n)$, the statement is equivalent to
$$\sin ((2j-1)\pi/n) >  2 \sin^2(j \pi/n) ,$$
 which,
using $2\sin^2(x)= 1- \cos(2x)$, is equivalent to
\begin{equation} \sin((2j-1)\pi/n) + \cos(2 j \pi/n) >1 \quad \text{for } 1 \le j \le \frac {n-1}4.
\label{e101}
\end{equation}
We will prove \eqref{e101} not only for integer $j$, but for all real $j \in [1, \frac {n-1}4]$.

The function $f(j)= \sin((2j-1)\pi/n) + \cos(2 j \pi/n)$ has the second derivative
$$ f''(j) = - (2\pi/n)^2 \left ( \sin((2j-1)\pi/n)  + \cos(2 j \pi/n) \right)$$
which is strictly negative on the interval $[1,\frac {n-1}4]$. Hence $f(j)$ achieves absolute minimum at one of the end points $1$ and $\frac {n-1}4$.
It is enough to show that the values of $f$ at these two end points are bigger than 1.

End point $1$. The inequality $f(1) > 1$ is

\begin{equation}
 \sin (\pi/n) + \cos (2 \pi/n)  >1. \label{e102}
 \end{equation}

The  function $f_1(x) = \sin x + \cos (2x)$ has the second derivative
$$ f''_1(x) = -\sin x - 4 \cos (2x)$$
which is strictly negative on the interval $(0,\pi/6)$. Hence on the closed interval $[0,\pi/6]$ the function $f_1(x)$ achieves the absolute minimum at one of the end points.
But $f_1(0)=f_1(\pi/6)=1$. If $n\ge 7$, then $\pi/n \in (0, \pi/6)$. Hence $f_1(\pi/n) >1$, which is \eqref{e102}.

End point $\frac {n-1}4$. One has 
\begin{eqnarray*} f\left (\frac {n-1}4\right)  & = &    \sin \left(\frac \pi 2 - \frac {3\pi}{2n}  \right) + \cos \left ( \frac \pi 2 - \frac \pi{2n}\right) \\
&=& \cos(3\pi/ 2n) + \sin (\pi/2n) \quad \text{ because } \sin(\pi /2 - x) = \cos x.
\end{eqnarray*}

Hence $f((n-1)/4) >1$ is equivalent to
\begin{equation}
 \sin (\pi/2n) + \cos (3 \pi/2n)  >1. \label{e102a}
 \end{equation}

Look at the function  $f_2(x) = \sin x + \cos (3x)$ on the interval $[0, \pi/14]$. The second derivative 
$$ f''_2(x) = -\sin x - 9 \cos (3x)$$
is strictly negative on the interval $(0,\pi/14)$, and $f_2(0)=1$ and $f_2(\pi/14)=1.004... >1$. It follows that $f_2(x) >1$ for $x \in (0, \pi/14]$.
If $n\ge 7$, then $\pi/2n \in (0, \pi/14]$. Hence $f_2(\pi/2n) >1$, which is \eqref{e102a}.
\end{proof}

\bibliographystyle{hamsalpha}\bibliography{biblio}
\end{document}